\documentclass[preprint,12pt]{elsarticle}

\usepackage{amsmath, amsthm, amssymb}
\usepackage{graphicx}
\usepackage{psfrag}
\usepackage{epsfig}
\usepackage{epstopdf}

\newtheorem{theorem}{Theorem}

\newcommand{{\Reals}}{\mathbb{R}}
\newcommand{\tn}{{\tilde{n}}}

\begin{document}

\begin{frontmatter}

\title{A Modular, Operator Splitting Scheme for Fluid-Structure Interaction Problems with Thick Structures}

 \author[martina]{M. Buka\v{c}}
 \author[UH]{S. \v{C}ani\'{c}}
 \author[UH] {R. Glowinski}
 \author[boris] {B. Muha}
 \author[UH]{A. Quaini}
 \address[martina]{Department of Mathematics, University of Pittsburgh, USA}
 \address[UH]{Department of Mathematics, University of Houston, USA}
 \address[boris]{Department of Mathematics, University of Zagreb, Croatia}




\begin{abstract}
We present an operator-splitting scheme for fluid-structure interaction (FSI) problems in hemodynamics,
where the thickness of the structural wall is comparable to the radius of the cylindrical fluid domain.
The equations of linear elasticity are used to model the structure, while the Navier-Stokes equations for
an incompressible viscous fluid are used to model the fluid. The operator splitting scheme,
based on Lie splitting, separates the elastodynamics structure problem,
from a fluid problem in which structure inertia is included to achieve unconditional stability. 
We prove energy estimates associated with unconditional stability of this modular scheme
for the full nonlinear FSI problem defined on a moving domain, without requiring any sub-iterations within time steps. 
Two numerical examples are presented, showing 
excellent agreement with the results of monolithic schemes. First-order convergence in time is shown numerically.
Modularity, unconditional stability without temporal sub-iterations, and simple implementation are the features that make this 
operator-splitting scheme particularly appealing for multi-physics problems involving fluid-structure interaction. 
\end{abstract}


\end{frontmatter}

\section{Introduction}
Fluid-structure interaction (FSI) problems arise in many applications.
They include  multi-physics problems in engineering such as aeroelasticity and propeller turbines, 
as well as biofluidic application such as self-propulsion organisms, fluid-cell interactions, and 
the interaction between blood flow and cardiovascular tissue.  
In biofluidic applications, such as the interaction between blood flow and cardiovascular tissue,
the density of the structure (arterial walls) is roughly equal to the density of the fluid (blood).
In such problems the energy exchange between the fluid and the structure is significant, leading to 
a highly nonlinear FSI coupling. 
A comprehensive study of these problems remains to be
a challenge due to their strong nonlinearity and multi-physics nature. 

The development of numerical solvers for fluid-structure interaction problems has become 
particularly active since the 1980's \cite{peskin1,peskin2,fauci1,fogelson,peskin3,peskin4,Griffith1,Griffith2,Griffith3,Griffith4,
donea1983arbitrary, hughes1981lagrangian,heil2004efficient,le2001fluid,leuprecht2002numerical,quarteroni2000computational,quaini2007semi,
baaijens2001fictitious,van2004combined,
fang2002lattice,feng2004immersed,krafczyk1998analysis,krafczyk2001two,
cottet2008eulerian,
figueroa2006coupled}.

\if 1 = 0
Among the most popular techniques are
the Immersed Boundary Method~\cite{peskin1,peskin2,fauci1,fogelson,peskin3,peskin4,Griffith1,Griffith2,Griffith3,Griffith4}
and the Arbitrary Lagrangian Eulerian (ALE) method
 \cite{donea1983arbitrary, hughes1981lagrangian,heil2004efficient,le2001fluid,leuprecht2002numerical,quarteroni2000computational,quaini2007semi}.
We further mention the Fictitious Domain Method in
combination with the mortar element method or ALE approach
\cite{baaijens2001fictitious,van2004combined}, and the methods recently proposed for the
use in the blood flow application such as
the Lattice Boltzmann method \cite{fang2002lattice,feng2004immersed,krafczyk1998analysis,krafczyk2001two},
the Level Set Method~\cite{cottet2008eulerian}
and the Coupled Momentum Method \cite{figueroa2006coupled}.
\fi

Until recently, only monolithic algorithms seemed  applicable to blood
flow simulations
\cite{figueroa2006coupled,gerbeau,nobile2008effective,zhao,bazilevs,bazilevs2}.
These algorithms are based on solving the entire nonlinear coupled problem as one monolithic system.
They are, however,
generally quite expensive in terms of computational time, programming time and
memory requirements, since they require solving a sequence of
strongly coupled problems using, e.g.,
the fixed point and Newton's methods
\cite{cervera,nobile2008effective,deparis,FernandezLifeV,heil2004efficient,matthies}.

The multi-physics nature of the blood flow problem
strongly suggest to employ partitioned (or staggered) numerical algorithms, where the coupled fluid-structure interaction
problem is separated into a fluid and a structure sub-problem.
The fluid and structure sub-problems are integrated in time in an alternating way, and the coupling conditions are enforced asynchronously. 
When the density of the structure is much larger than the density of the fluid, as is the case in aeroelasticity, 
it is sufficient to solve, at every time step, just one fluid sub-problem and one structure sub-problem to obtain a solution.
The classical {\ loosely-coupled} partitioned schemes of this kind typically use the structure velocity in the {\ fluid sub-problem} as 
{Dirichlet} data for the fluid velocity (enforcing the no-slip boundary condition at the fluid-structure interface),
while in the { structure sub-problem} the structure is loaded by the fluid {normal stress} calculated in the fluid sub-problem.
These {{Dirichlet-Neumann}} loosely-coupled partitioned schemes work well for problems in which the structure is much heavier than the fluid. 
Unfortunately, when fluid and structure have comparable densities, which is  the case in the blood flow application, the simple strategy of separating the fluid from the structure suffers from severe stability 
issues \cite{causin2005added,Michler}. This is because the energy of the discretized problem in Dirichlet-Neumann loosely-coupled schemes
does not approximate well the energy of the continuous problem. 
A partial solution to this problem is to sub-iterate several times between the fluid and structure sub-solvers
at every time step until the energy of the continuous problem is well approximated.
These {\ strongly-coupled} partitioned schemes, however,  are computationally expensive and may suffer from
convergence issues for certain parameter values \cite{causin2005added}.

To get around these difficulties, and to retain the main advantages of loosely-coupled partitioned schemes
such as modularity, simple implementation, and low computational costs, several new loosely-coupled algorithms 
have been proposed recently. In general, they behave quite well for FSI problems containing a thin fluid-structure interface with mass
\cite{badia,Martina_paper1,Martina_Multilayered,MarSunLongitudinal,guidoboni2009stable,nobile2008effective,Fernandez1,Fernandez2,Fernandez3,Fernandez2006projection,astorino2009added,astorino2009robin,badia2008splitting,quaini2007semi,murea2009fast,deparis,deparis2}.

For FSI problems in which the structure is ``thick'', i.e., the thickness of the structure is comparable to the 
transverse dimension of the fluid domain, partitioned, loosely-coupled schemes
 are more difficult to construct. In fact, to the best of our knowledge, there have been no loosely-coupled, partitioned schemes
 proposed so far in literature for a class of FSI problems in hemodynamics that contain thick structure models to study the 
 elastodynamics of arterial walls. The closest works in this direction include a strongly-coupled partition scheme by Badia et al. in \cite{badia2009robin},
 and an explicit scheme by Burman and Fern\'{a}ndez where certain ``defect-correction'' sub-iterations are necessary to achieve optimal accuracy
 \cite{burman2009stabilization}. 

 More precisely, in \cite{badia2009robin}, the authors construct a strongly-coupled partitioned scheme based on certain Robin-type coupling conditions. In addition to the classical Dirichlet-Neumann and Neumann-Dirichlet schemes, they
also propose a Robin-Neumann and a Robin-Robin scheme, that converge without relaxation, and
need a smaller number of sub-iteration between the fluid and the
structure in each time step than classical strongly-coupled schemes. 

In \cite{burman2009stabilization}, Burman and Fern\'{a}ndez propose an explicit scheme where 
the coupling between the fluid and a thick structure is enforced in a weak sense using Nitsche's approach \cite{hansbo2005nitsche}.
The formulation in \cite{burman2009stabilization} still suffers from stability issues related to the added mass effect, which were
corrected by adding a weakly consistent penalty term that includes pressure variations at the interface. The added term, however,
lowers the temporal accuracy of the scheme, which was then corrected by proposing a few defect-correction sub-iterations to achieve
optimal accuracy.


In the work presented here, we take a different approach to separate the calculation of the fluid and structure
sub-problems in the case when the FSI problem incorporates a thick elastic structure. Our approach is based on Lie splitting,
also known as the Marchuk-Yanenko scheme. This splitting is applied to the coupled FSI problem written in ALE form. Namely,
to deal with the motion of the fluid domain, in this manuscript we utilize an Arbitrary Lagrangian Eulerian (ALE) approach
\cite{hughes1981lagrangian,donea1983arbitrary,heil2004efficient,leuprecht2002numerical,quaini2007semi,quarteroni2000computational,le2001fluid}. 
Once the coupled problem is written in ALE form, the Lie splitting is applied.
The coupled FSI problem in ALE form  is split into a fluid sub-problem, and a structure sub-problem.
The fluid sub-problem includes structure inertia
 to avoid instabilities associated with the added mass effect in partitioned schemes.
 This also avoids the need for any sub-iterations in each time step.
 A structure elastodynamics problem is then  solved separately. 
We first introduced this approach in \cite{guidoboni2009stable} to deal with FSI problems containing thin structures, leading to
a completely partitioned, loosely-coupled scheme called the kinematically-coupled scheme. To increase the accuracy of the kinematically-coupled scheme, 
in \cite{Martina_paper1} we introduced a modification of this scheme, called the kinematically-coupled $\beta$-scheme,
which was based on including a $\beta$-fraction of the pressure at the fluid-structure interface into the structure sub-problem. 
Another novelty of \cite{Martina_paper1} was the fact that 
this scheme was applied to a FSI interaction problem where the structure was modeled by the Koiter shell model accounting for both
radial and longitudinal displacements.  Due to its simple implementation, modularity, and good performance,
modifications of this scheme have been used by several authors to study FSI problems
in hemodynamics including 
cardiovascular stents \cite{BorSunStent},  thin structures with longitudinal displacement \cite{MarSunLongitudinal},
multi-layered structure of arterial walls \cite{BorSunMulti,Martina_Multilayered}, 
poroelastic arterial walls \cite{Martina_Biot}, or non-Newtonian fluids \cite{Lukacova}.
In the present paper we extend the kinematically-coupled $\beta$-scheme to FSI problems with thick structures.

This extension is not trivial because the resulting scheme, unlike those cited above, is not completely partitioned due to the fact
that in problems with thick structures, the {\sl fluid-structure interface does not have a well-defined mass/inertia}. 
More precisely, to achieve unconditional stability, our operator splitting strategy is based on 
including the fluid-structure interface inertia into the fluid sub-problem. 
This can be easily done when the fluid-structure interface has mass. In that case the structure inertia can be included in the fluid sub-problem
through a Robin-type boundary condition at the fluid-structure interface \cite{Martina_paper1}.
However, in problems in which the interface between the fluid and structure is just a trace of a thick structure that is in contact with the fluid,
as is the case in the present manuscript, the inclusion of the structure inertia in the fluid
sub-problem is problematic if one wants to split the problem in the spirit of partitioned schemes.
We address this issue by defining a new ``fluid sub-problem'' which involves solving a {\sl simplified coupled} problem on both the fluid domain and the structure domain, {\sl in a monolithic fashion}. The inertia of the structure is included in this ``fluid sub-problem'' 
not through a boundary condition for the fluid sub-problem, but by solving a simple, structure problem
involving only structure inertia (and structural viscosity if the structure is {\sl visco}elastic), coupled 
with the fluid problem via a simple continuity of stress condition at the interface.
Although solving this simplified coupled problem on both domains
is reminiscent of monolithic FSI schemes, the situation is, however, {\sl much simpler}, since
the {\sl hyperbolic effects associated with fast wave propagation 
in the structural elastodynamics problem are not included here}. As a result, we show below in
Section~\ref{condition_number}, that the condition number of this sub-problem
is smaller by several orders of magnitude than the condition number associated with monolithic FSI schemes.
In fact, the condition number of this  sub-problem is of the same order of magnitude 
as the condition number of the pure fluid sub-problem. 
Furthermore,  the time step in this sub-problem can be taken larger than the time step in the classical monolithic schemes,
which is dictated by the fast traveling waves in the elastic structure. 
Using this approach we achieved unconditional stability of this operator splitting scheme that separates the fluid from the structure sub-problems
without a need for sub-iterations, but, as mentioned above, the drawback is the expense of generating the computational mesh on both, the fluid and structure domains to resolve the fluid sub-problem. 

In Section~\ref{stability} we prove that for the fully nonlinear FSI problem defined on a moving domain, the proposed scheme
satisfies an energy estimate which is associated with unconditional stability of the scheme, for all the parameters in the problem. 
To the best of our knowledge,
this is the first result in which an energy estimate indicating unconditional stability of a partition-like scheme
is obtained  for a full nonlinear FSI problem.

In Section~\ref{sec5} we study two examples from literature involving FSI problems in hemodynamics with thick structures.
We showed that in both cases our simulations compared well with the results of monolithic schemes. 
Furthermore, in Section~\ref{accuracy} we showed, on a numerical example, that our scheme is first-order accurate in time. 

Although the presentation and numerical examples in this manuscript are given in terms of 2D problems, there is nothing in the 
operator-splitting scheme 
that depends on the dimension of the problem. Therefore, the same ideas as those presented in this manuscript
apply to problems in 3D.

We conclude by emphasizing that, as in partitioned schemes, our scheme is modular in the sense that 
different modules can be easily added or replaced to study more complex multi-physics problems,
and no sub-iterations between the different modules are needed to achieve stability. 
Thus, we argue that the proposed operator splitting scheme is closer in spirit to partitioned (loosely-coupled) schemes than to monolithic schemes, 
providing an appealing approach to solve coupled FSI in hemodynamics with thick structures.

\section{Description of the fluid-structure interaction problem}\label{sec3}

\begin{figure}[ht]
\centering{
\includegraphics[scale=0.57]{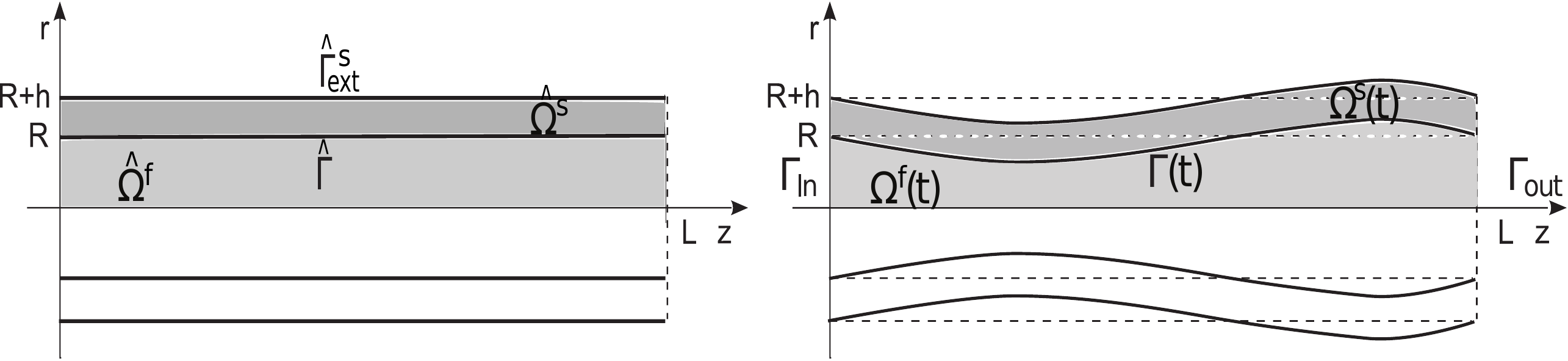}
}
\caption{Left: Reference domains $\hat{\Omega}^f \cup \hat{\Omega}^s$. Right: Deformed domains $\Omega^f(t) \cup \Omega^s(t).$}
\label{fig:domain}
\end{figure}

We consider the flow of an incompressible, viscous fluid in a two-dimensional channel of reference length $L$, and reference width $2R$,
see Figure~\ref{fig:domain}.
The lateral boundary of the channel is bounded by a thick, deformable wall with finite thickness $h$,
i.e., the wall thickness $h$ is not necessarily small with respect to the channel radius $R$.
We are interested in simulating a pressure-driven flow through a deformable 2D channel, 
in which the fluid and structure are fully coupled via a set of coupling conditions describing two-way coupling.
Without loss of generality, we consider only the upper half of the fluid domain supplemented by a symmetry condition at the axis of symmetry.
Thus, the reference fluid and structure domains in our problem are given, respectively, by
\begin{eqnarray*}
\hat{\Omega}^f &:=& \{(z,r) | 0<z<L, 0<r<R \}, \\
\hat{\Omega}^s &:=& \{(z,r) | 0<z<L, R<r<R+h \},
\end{eqnarray*}
with $\hat\Gamma$ denoting the Lagrangian boundary of $\hat\Omega^s$ in contact with the fluid $$\hat\Gamma=(0,L).$$
Here $z$ and $r$ denote the horizontal and vertical Cartesian coordinates, respectively (see Figure~\ref{fig:domain}). 
Throughout the rest of the manuscript we will be using the ``hat'' notation to denote the Lagrangian variables defined on the reference configuration.

The flow of an incompressible, viscous fluid is modeled by  the Navier-Stokes equations:
\begin{eqnarray}\label{NS1}
 \rho_f \bigg( \frac{\partial \boldsymbol{v}}{\partial t}+ \boldsymbol{v} \cdot \nabla \boldsymbol{v} \bigg) & =& \nabla \cdot \boldsymbol\sigma \quad \;  \textrm{in}\; \Omega^f(t)\; \textrm{for}\; t \; \in(0,T), \\
 \label{NS2}
\nabla \cdot \boldsymbol{v} &=& 0 \quad \quad \; \quad \textrm{in}\; \Omega^f(t)\; \textrm{for}\; t \; \in(0,T),
\end{eqnarray}
where $\boldsymbol{v}=(v_z,v_r)$ is the fluid velocity, $p$ is the fluid pressure, $\rho_f$ is the fluid density, and $\boldsymbol\sigma$ is the fluid Cauchy stress tensor. 
For a Newtonian fluid the Cauchy stress tensor is given by $\boldsymbol\sigma = -p \boldsymbol{I} + 2 \mu_f \boldsymbol{D}(\boldsymbol{v}),$ where
$\mu_f$ is the fluid viscosity and  $\boldsymbol{D}(\boldsymbol{v}) = (\nabla \boldsymbol{v}+(\nabla \boldsymbol{v})^{\tau})/2$ is the rate-of-strain tensor.

The fluid domain $\Omega^f(t)$ is not known {\sl a priori} as it depends on the solution of the problem. 
Namely, the lateral boundary of $\Omega^f(t)$ is determined by the trace of the displacement of the thick structure
at the fluid-structure interface, as it will be shown below.

At the inlet and outlet boundaries of the fluid domain, denoted by  $\Gamma^f_{in}$ and $\Gamma^f_{out},$ respectively,
we prescribe the normal stress:
\begin{eqnarray}
\boldsymbol {\sigma n}^f_{in}(0,r,t) &=& -p_{in}(t) \boldsymbol{n}^f_{in} \quad \textrm{on} \; (0,R) \times (0,T), \label{inlet} \\
\boldsymbol {\sigma n}^f_{out}(L,r,t) &=& -p_{out}(t) \boldsymbol{n}^f_{out} \quad \textrm{on} \; (0,R) \times (0,T),\label{outlet}
\end{eqnarray}
where $\boldsymbol{n}^f_{in}/\boldsymbol{n}^f_{out}$ are the outward normals to the inlet/outlet fluid boundaries, respectively.
Even though not physiologically optimal, these boundary conditions are common in blood flow modeling~\cite{badia2008fluid, miller2005computational, nobile2001numerical}.

At the bottom fluid boundary $r=0$, denoted by $\Gamma_b$, the following symmetry conditions are prescribed:
\begin{equation}\label{symmetry_condition}
 \frac{\partial v_z}{\partial r}(z,0,t) = 0, \quad v_r(z,0,t) = 0 \quad \textrm{on} \; (0,L) \times (0,T).
\end{equation}

The lateral fluid boundary is bounded by a deformable, two-dimensional, 
elastic/viscoelastic  wall with finite thickness, modeled by the following elastodynamics equations:
\begin{equation}
\rho_s \frac{\partial^2 \hat{\boldsymbol U}}{\partial t^2} + \gamma \hat{\boldsymbol U}= \nabla \cdot \boldsymbol S( \hat{\boldsymbol U}) + \epsilon \Delta \frac{\partial \hat{\boldsymbol U}}{\partial t}\quad \textrm{in} \; \hat{\Omega}^s \; \textrm{for} \; t \in (0,T),
\end{equation}
where $\hat{\boldsymbol U}(\hat{z},\hat{r},t) = (\hat{U}_z(\hat{z},\hat{r},t), \hat{U}_r(\hat{z},\hat{r},t))$ is the structure displacement, $\rho_s$ is the structure density, and $\epsilon$ is the constant
modeling viscoelastic structural effects. 
In the work described here, $\epsilon$ can be taken to be zero, thereby accounting for 
a possibility of modeling strictly elastic structures. 
The term $\gamma {\hat{\boldsymbol U}}$ (i.e., the linearly elastic spring term)
 comes from 3D axial symmetry, accounting for the recoil due to circumferential strain, keeping the top and bottom
structure displacement connected in 2D, see, e.g., \cite{badia2008splitting,ma1992numerical,barker2010scalable}. 
Tensor ${\bf S}$ is the first Piola-Kirchhoff stress tensor given by $\boldsymbol S( \hat{\boldsymbol U}) = 2 \mu_s \boldsymbol D(\hat{\boldsymbol U}) + \lambda_s (\nabla \cdot \hat{\boldsymbol U}) \boldsymbol I$, where $\mu_s$ and $\lambda_s$ are the Lam\'e constants. 
The structure is described in the Lagrangian framework, defined on a fixed, reference domain $\hat\Omega^s$.
In contrast, the fluid problem, written in the Eulerian framework, is defined on a domain $\Omega^f(t)$ which depends on time.

We assume that the structure is fixed at the inlet and outlet portion of the boundary:
\begin{equation}\label{homostructure1}
 \hat{\boldsymbol U}(0, \hat{r}, t) = \hat{\boldsymbol U}(L, \hat{r}, t) = 0, \quad \textrm{for} \; \hat{r} \in [R, R+h], t \in (0,T).
\end{equation}
At the external structure boundary $\hat{\Gamma}^s_{ext}$, we assume that
the structure is exposed to a certain external ambient pressure $P_{ext}$, while the axial displacement remains fixed:
\begin{eqnarray}\label{homostructure2}
 \boldsymbol n^s_{ext} \cdot \boldsymbol {S n}^s_{ext} &=&  -P_{ext}, \quad \textrm{on} \; \hat{\Gamma}^s_{ext} \times (0,T), \\
 \hat{U}_z &=&  0, \quad \quad \quad \textrm{on} \; \hat{\Gamma}^s_{ext} \times (0,T),
\end{eqnarray}
where $\boldsymbol n^s_{ext}$ is the outward unit normal vector on $\hat{\Gamma}^s_{ext}$.

Initially, the fluid and the structure are assumed to be at rest, with zero displacement from the reference configuration
\begin{equation}\label{initial}
 \boldsymbol{v}=0, \quad \hat{\boldsymbol U} = 0, \quad \frac{\partial \hat{\boldsymbol U}}{\partial t}=0.
 \end{equation}

The fluid and structure are coupled via the kinematic and dynamic boundary conditions~\cite{canic2005self,mikelic2007fluid}:
\begin{itemize}
 \item \textbf{Kinematic coupling condition} 
describes continuity of velocity
\begin{equation}\label{kinematic}
 \boldsymbol{v}(\hat{z}+\hat{U}_z(\hat{z}, R,t),R+\hat{U}_r(\hat{z}, R,t),t)=\frac{\partial \hat{\boldsymbol U}}{\partial t}(\hat{z}, R, t) \quad \; \textrm{on} \; (0,L)\times (0,T).
\end{equation} 
\item \textbf{Dynamic coupling condition} 
describes balance of contact forces:
\begin{equation}
 J\ \widehat{\boldsymbol{\sigma n}^f}|_{\hat\Gamma} +  \boldsymbol {Sn}^s|_{\hat{\Gamma}} +\epsilon \frac{\partial}{\partial \boldsymbol n^s} \bigg(\frac{\partial \hat{\boldsymbol U}}{\partial t} \bigg) \bigg|_{\hat{\Gamma}}= 0 \quad \textrm{on} \; (0,L)\times(0,T),   
\label{dynamic}
\end{equation}
where $J$ denotes the Jacobian of the transformation from the Eulerian to Lagrangian framework, 
$J = \sqrt{(1+\partial \eta_z/\partial z)^2 + (\partial \eta_r/\partial z)^2}$,
and $\widehat{\boldsymbol{\sigma n}^f}$ denotes the normal fluid stress defined on $\hat{\Omega}^f = (0,L)\times(0,R)$ (here ${\boldsymbol{n}}^f$ is the outward unit normal
to the deformed fluid domain), and $\boldsymbol n^s$ is the outward unit normal to the structural domain.
\end{itemize}

\subsection{The ALE framework}
\begin{figure}[ht]
 \centering{
 \includegraphics[scale=0.35]{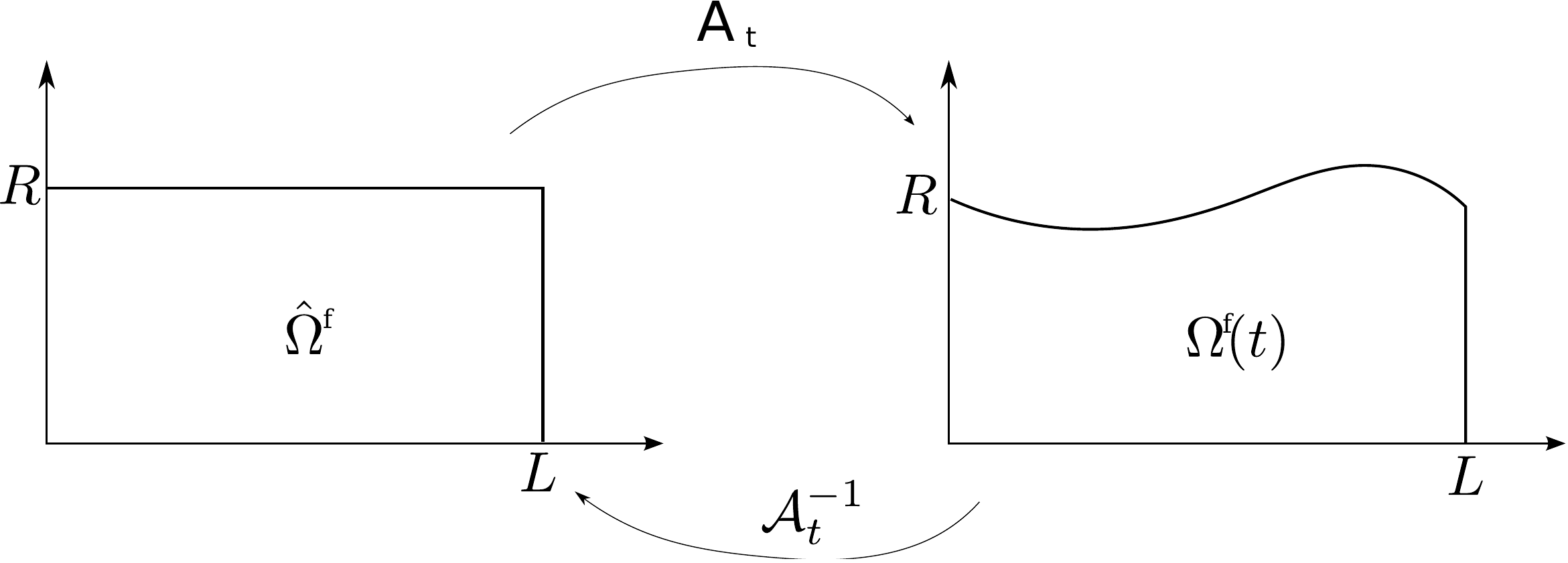}
 }
 \caption{$\mathcal{A}_t$ maps the reference domain $\hat\Omega^f$ into the current domain $\Omega^f(t)$. }
\label{ale703}
 \end{figure}

To deal with the motion of the fluid domain we adopt the Arbitrary Lagrangian-Eulerian (ALE) approach~\cite{hughes1981lagrangian,donea1983arbitrary,nobile2001numerical}.
In the context of finite element method approximation of moving-boundary problems,
ALE method deals efficiently 
with the deformation of the mesh, especially near the interface
between the fluid and the structure,
and with the issues 
related to the approximation of the time-derivatives
$\partial \boldsymbol v/\partial t \approx (\boldsymbol {v}(t^{n+1})-\boldsymbol {v}(t^{n}))/\Delta t$ which,
due to the fact that $\Omega^f(t)$ depends on time, is not well defined since
the values $\boldsymbol {v}(t^{n+1})$ and $\boldsymbol {v}(t^{n})$ 
correspond to the values of $\boldsymbol {v}$ defined at two different domains.
ALE approach is based on introducing
a family of (arbitrary, invertible, smooth) mappings ${\cal{A}}_t$ defined on 
a single, fixed, reference domain $\hat{\Omega}^f$ such that, for each 
$t \in (t_0, T)$, ${\cal{A}}_t$ maps the reference domain $\hat{\Omega}^f=(0,L) \times (0,R)$ into the current domain $\Omega^f(t)$ (see Figure~\ref{ale703}):
$$\mathcal{A}_t : \hat{\Omega}^f \subset \mathbb{R}^2 \rightarrow \Omega^f(t) \subset \mathbb{R}^2,  \quad \boldsymbol{x}=\mathcal{A}_t(\hat{\boldsymbol{x}}) \in \Omega^f(t), \quad \textrm{for} \; \hat{\boldsymbol{x}} \in \hat{\Omega}^f.$$
In our approach, we define $\mathcal{A}_t$ to be a harmonic extension of the structure displacement $\hat{\boldsymbol U}$ at the fluid-structure interface onto the whole domain $\hat\Omega^f$, for a given $t:$
\begin{eqnarray*}
 \Delta \mathcal{A}_t &=& 0 \quad \rm{in} \; \hat{\Omega}^f, \\
\mathcal{A}_t |_{\hat{\Gamma}} &=& \hat{\boldsymbol U}|_{\hat{\Gamma}}, \\
\mathcal{A}_t |_{\partial \hat{\Omega}^f\backslash \hat{\Gamma}}&=&0.
\end{eqnarray*}

To rewrite system~\eqref{NS1}-\eqref{NS2} in the ALE framework
we notice that for a function $f=f(\boldsymbol {x},t)$ defined on $\Omega^f(t) \times (0,T)$ 
the corresponding function $\hat{f} := f \circ \mathcal{A}_t$ defined on $\hat{\Omega} \times (0,T)$
is given by
$$\hat{f}(\hat{\boldsymbol{x}},t) = f(\mathcal{A}_t(\hat{\boldsymbol{x}}),t).$$
Differentiation with respect to time, after using the chain rule, gives
\begin{equation}
 \frac{\partial f}{\partial t}\bigg|_{\hat{\boldsymbol x}} =  \frac{\partial {f}}{\partial t}+\boldsymbol{w} \cdot \nabla f,
\end{equation}
where 
$\boldsymbol{w}$ denotes domain velocity given by
\begin{equation}
 \boldsymbol{w}(\boldsymbol {x},t) = \frac{\partial \mathcal{A}_t(\hat{\boldsymbol x})}{\partial t}.\bigg. \label{w}
\end{equation}

Finally, system~\eqref{NS1}-\eqref{NS2} in ALE framework read as follows: 
Find $\boldsymbol{v}= (v_z,v_r)$ and $p$,  
with $\hat{\boldsymbol{v}}(\hat{\boldsymbol{x}},t) = \boldsymbol{v}(\mathcal{A}_t(\hat{\boldsymbol{x}}),t)$  such that
\begin{eqnarray}
& & \rho_f \bigg( \frac{\partial \boldsymbol{v}}{\partial t}\bigg|_{\hat{\boldsymbol x}}+ (\boldsymbol{{v}}-\boldsymbol{w}) \cdot \nabla \boldsymbol{v} \bigg)  = \nabla \cdot \boldsymbol\sigma,   \quad \quad \textrm{in}\; \Omega^f(t) \times (0,T), \label{ALEsys1}\\
& &\nabla \cdot \boldsymbol{v} = 0  \quad \quad \; \textrm{in}\; \Omega^f(t) \times (0,T),
\label{ALEsys2}
\end{eqnarray}
with corresponding initial and boundary conditions.

We will apply the Lie splitting on this system written in ALE form. Before we do that, 
we introduce the notion of weak solutions to the coupled FSI problem. For this purpose, 
we utilize the ALE mapping, introduced above, and define the function spaces on moving domains
in terms of the ALE-mapped functions, defined on the fixed, reference domain.

\subsection{Weak solution of the coupled FSI problem}
Define the following test function spaces  
\begin{eqnarray}
V^f(t) &=& \{\boldsymbol {\varphi}: \Omega^f(t) \rightarrow \mathbb{R}^2| \; \boldsymbol \varphi = \hat{\boldsymbol \varphi} \circ (\mathcal{A}_t)^{-1}, \hat{\boldsymbol \varphi} \in (H^1(\hat{\Omega}^f))^2,  \nonumber \\
& & \quad \quad  \varphi_r|_{r=0}=0,  \}, \label{Vf(t)} \\
Q(t) &=& \{q: \Omega^f(t) \rightarrow \mathbb{R}| \; q = \hat{q} \circ (\mathcal{A}_t)^{-1}, \hat{q} \in L^2(\hat{\Omega}^f) \}, \label{Q(t)} \\
 \hat{V}^s &=& \{\hat{\boldsymbol \psi}: \hat{\Omega}^s \rightarrow \mathbb{R}^2| \; \hat{ \boldsymbol \psi} \in (H^1(\hat{\Omega}^s))^2, \hat{ \boldsymbol \psi} |_{z=0, L}=0, \hat{\psi}_z|_{\hat{\Gamma}^s_{ext}} = 0  \}, \label{Vs(t)}
\end{eqnarray}
for all $t \in [0,T)$, and introduce the following function space 
\begin{equation}\label{Vfsi}
V^{fsi}(t) = \{(\boldsymbol \varphi, \hat{ \boldsymbol \psi}) \in V^f(t) \times \hat{V}^s| \; \boldsymbol \varphi|_{\Gamma(t)} = \hat{ \boldsymbol \psi}|_{\hat{\Gamma}} \circ (\mathcal{A}_t^{-1}|_{\Gamma(t)})\}.
\end{equation}
The variational formulation of the coupled fluid-structure interaction problem now reads: 
For $t \in (0, T)$, find $(\boldsymbol v, \hat{\boldsymbol U}, p) \in V^{f}(t) \times \hat{V}^s \times Q(t)$ such that ${\boldsymbol v = \hat{\boldsymbol U}}_t\circ \mathcal{A}_t^{-1}$ on $\Gamma(t)$,
and the following holds for all $(\boldsymbol \varphi, \hat{ \boldsymbol \psi}, q) \in V^{fsi}(t) \times Q(t)$:
\begin{equation*}
\rho_f \left\{\int_{\Omega^f(t)} \frac{\partial \boldsymbol v}{\partial t} \cdot \boldsymbol \varphi d\boldsymbol x+ \int_{\Omega^f(t)} (\boldsymbol v \cdot \nabla) \boldsymbol v \cdot \boldsymbol \varphi d \boldsymbol x \right\} +2 \mu_f \int_{\Omega^f(t)} \boldsymbol D(\boldsymbol v) : \boldsymbol D(\boldsymbol \varphi) d \boldsymbol x 
\end{equation*}
\begin{equation*}
- \int_{\Omega^f(t)} p \nabla \cdot \boldsymbol \varphi d \boldsymbol x +\rho_s \int_{\hat{\Omega}^s} \frac{\partial^2 \hat{\boldsymbol U}}{\partial t^2} \cdot \hat{ \boldsymbol \psi} d\boldsymbol x+2 \mu_s \int_{\hat{\Omega}^s} \boldsymbol D(\hat{\boldsymbol U}) : \boldsymbol D(\hat{ \boldsymbol \psi}) d \boldsymbol x 
\end{equation*}
\begin{equation*}
+\lambda_s \int_{\hat{\Omega}^s} (\nabla \cdot \hat{\boldsymbol U}) (\nabla \cdot \hat{ \boldsymbol \psi}) d \boldsymbol x +\gamma \int_{\hat{\Omega}^s} \hat{\boldsymbol U} \cdot \hat{ \boldsymbol \psi} d\boldsymbol x+ \epsilon \int_{\hat{\Omega}^s} \nabla \frac{\partial \hat{\boldsymbol U}}{\partial t} : \nabla \hat{ \boldsymbol \psi} d\boldsymbol x
\end{equation*}
\begin{equation*}
 + \int_{\hat{\Gamma}^s_{ext}} P_{ext} (\hat{ \boldsymbol \psi} \cdot \boldsymbol n^s_{ext}) dS= \int_0^R p_{in}(t) \varphi_z|_{z=0} dr-  \int_0^R p_{out}(t) \varphi_z|_{z=L} dr,
\end{equation*}
with 
\begin{equation*}
\int_{\Omega^f(t)} q \nabla \cdot \boldsymbol v d \boldsymbol x = 0.
\end{equation*}

\section{The Numerical Scheme}\label{sec4}

\subsection{An Operator Splitting Approach}
To solve the fluid-structure interaction problem written in ALE form~\eqref{ALEsys1}-\eqref{ALEsys2}, and
with the initial and boundary conditions
 \eqref{inlet}-\eqref{dynamic}, we propose an operator splitting scheme which is based on 
 Lie operator splitting~\cite{glowinski2003finite}, also known as 
Marchuk-Yanenko splitting.
The splitting is introduced to separate the different physics in the problem represented by the
fluid and structure sub-problems.
To achieve unconditional stability of such an approach, 
it is crucial to include inertia of the fluid-structure interface in the fluid sub-problem,
thereby avoiding the ``added mass effect'' responsible for 
instabilities of classical Dirichlet-Neumann partitioned schemes \cite{causin2005added}.
This can be easily done when the fluid-structure interface has mass, e.g., the interface is modeled
as a thin structure satisfying the Koiter shell equations, as was done in \cite{guidoboni2009stable,Martina_paper1}.
In that case the structure inertia can be included in the fluid sub-problem
through a Robin-type boundary condition at the fluid-structure interface \cite{Martina_paper1}.
However, in problems in which the interface between the fluid and structure is just a trace of a thick structure that is in contact with the fluid,
as is the case in the present manuscript, the inclusion of structure inertia in the fluid
sub-problem is problematic if one wants to split the problem in the spirit of partitioned schemes.
We address this issue by defining a new ``fluid sub-problem'' which involves solving a {\sl simplified coupled} problem on both the fluid domain and the structure domain. 
The fluid sub-problem is solved {\sl coupled with structure inertia}
(and with the viscosity of the structure if the structure is viscoelastic) in a monolithic way,
leaving out the elastic part of the structure. Although solving this simplified coupled problem on both domains
is reminiscent of monolithic FSI schemes, the situation is, however, much simpler, since
the hyperbolic effects associated with fast wave propagation 
in the structural elastodynamics problem are not included here. As a result, we show below in
Section~\ref{condition_number}, that the condition number of this sub-problem
is smaller by several orders of magnitude than the condition number associated with monolithic FSI schemes.
In fact, the condition number of this  sub-problem is of the same order of magnitude 
as the condition number of a pure fluid sub-problem!
This approach also allows the choice of larger time-steps, not restricted by the time-scale associated with 
fast wave propagation in the elastic structure. 

\subsubsection{Lie splitting and the first-order system}\label{LieSplittingGeneral}
The Lie operator splitting is defined for evolutionary problems which can be written as a first-order system in time:
\begin{eqnarray}\label{LieProblem}
   \frac{\partial \phi}{\partial t} + A(\phi) &=& 0, \quad \textrm{in} \ (0,T), \\
\phi(0) &=& \phi_0, 
\end{eqnarray}
where A is an operator from a Hilbert space into itself, and $A$ can be split, in a non-trivial decomposition, as
\begin{equation}
 A = \sum\limits_{i=1}^I A_i.
\end{equation}
The Lie scheme consists of the following. Let $\Delta t>0$ be a time discretization step. Denote $t^n=n\Delta t$ and let $\phi^n$
be an approximation of $\phi(t^n).$ Set $\phi^0=\phi_0.$ Then, for $n \geq 0$ compute $\phi^{n+1}$ by solving
\begin{eqnarray}
   \frac{\partial \phi_i}{\partial t} + A_i(\phi_i) &=& 0 \quad \textrm{in} \; (t^n, t^{n+1}), \\
\phi_i(t^n) &=& \phi^{n+(i-1)/I}, 
\end{eqnarray}
and then set $\phi^{n+i/I} = \phi_i(t^{n+1}),$ for $i=1, \dots. I.$
This method is first-order accurate in time. 
More precisely, if~\eqref{LieProblem} is defined on a finite-dimensional space, and if the operators $A_i$ are smooth enough, 
then $\| \phi(t^n)-\phi^n \| = O(\Delta t)$~\cite{glowinski2003finite}.

To apply the Lie operator splitting scheme, we must rewrite system~\eqref{NS1}-\eqref{dynamic} in first-order form. 
For this purpose we introduce the structure velocity $\hat{\boldsymbol V}=\frac{\partial \hat{\boldsymbol U}}{\partial t}$,
and rewrite the problem as follows:
Find $\boldsymbol{v}= (v_z,v_r)$, $p$,  $\hat{\boldsymbol U} = (\hat{U}_z, \hat{U}_r)$ and $\hat{\boldsymbol V} = (\hat{V}_z, \hat{V}_r)$,
with $\hat{\boldsymbol{v}}(\hat{\boldsymbol{x}},t) = \boldsymbol{v}(\mathcal{A}_t(\hat{\boldsymbol{x}}),t),$  such that
\begin{eqnarray}
 \rho_f \bigg(\displaystyle{ \frac{\partial \boldsymbol{v}}{\partial t}}\bigg|_{\hat{\boldsymbol x}}+ (\boldsymbol{{v}}-\boldsymbol{w}) \cdot \nabla \boldsymbol{v} \bigg)  &=& \nabla \cdot \boldsymbol\sigma,   \quad  \textrm{in}\; \Omega^f(t) \times (0,T), \label{sys1}\\ 
 \nabla \cdot \boldsymbol{v} &=& 0,  \  \; \textrm{in}\; \Omega^f(t) \times (0,T), \\
\rho_s \frac{\partial \hat{\boldsymbol V}}{\partial t} + \gamma \hat{\boldsymbol U} &=& \nabla \cdot \boldsymbol S( \hat{\boldsymbol U}) + \epsilon \Delta \hat{\boldsymbol V},
\   \textrm{in} \; \hat{\Omega}_s\times (0,T), \\
\frac{\partial \hat{\boldsymbol U}}{\partial t} &=& \hat{\boldsymbol V},  \   \textrm{in} \; \hat{\Omega}_s \times (0,T),
\end{eqnarray}
with the coupling conditions at the fluid-structure interface
\begin{eqnarray}
& & \hat{\boldsymbol V}|_{\hat{\Gamma}} = \hat{\boldsymbol{v}}|_{\hat\Gamma}, \\
& &  J\ \widehat{\boldsymbol{\sigma n}^f}|_{\hat\Gamma} +  \boldsymbol {Sn}^s|_{\hat{\Gamma}} +\epsilon \frac{\partial \hat{\boldsymbol V}}{\partial \boldsymbol n^s}\bigg|_{\hat{\Gamma}}= 0,
\end{eqnarray}
(where $\hat{\boldsymbol{v}}|_{\hat\Gamma}$ denotes  $\boldsymbol {v}(\hat{z} + \hat{U}_z(\hat{z},R,t), R + \hat{U}_r(\hat{z},R,t),t)$; 
and $\widehat{\boldsymbol{\sigma n}^f}|_{\hat\Gamma}$ denotes the normal fluid stress defined on the reference fluid domain $\hat\Omega^f$),
and with boundary conditions
\begin{eqnarray}
& & \frac{\partial v_z}{\partial r}(z,0,t) =  v_r(z,0,t) = 0 \quad \textrm{on} \; \Gamma^f_b, \\
& & \boldsymbol{v}(0,R,t) = \boldsymbol{v}(L,R,t) = 0, \quad \hat{\boldsymbol U}|_{\hat{z}=0,L}=0 , \\
& & \boldsymbol{\sigma n}^f_{in}(0, r,t) = -p_{in}(t) \boldsymbol{n}^f_{in}, \\
& & \boldsymbol{\sigma n}^f_{out}(L,r,t) = -p_{out}(t) \boldsymbol{n}^f_{out} \; \textrm{on} \; (0,R) \times (0,T). \\
& &  \boldsymbol n^s_{ext} \cdot \boldsymbol {S n}^s_{ext} =  -P_{ext} \quad \textrm{on} \; \hat{\Gamma}^s_{ext} \times (0,T), \\
& & \hat{U}_z =  0 \quad \quad \quad \textrm{on} \; \hat{\Gamma}^s_{ext} \times (0,T).
\end{eqnarray}
At time $t = 0$ the following initial conditions are prescribed:
\begin{eqnarray}
& & \boldsymbol{v}|_{t=0} = \boldsymbol{0}, \quad \hat{\boldsymbol U}|_{t=0} = \boldsymbol{0}, \quad \hat{\boldsymbol V} |_{t=0} = \boldsymbol{0}.  \label{sys2}
\end{eqnarray}
Using the Lie operator splitting approach,
problem \eqref{sys1}-\eqref{sys2} is split into a sum of the following sub-problems: 
\begin{enumerate}
 \item[{\bf A1.}] An elastodynamics sub-problem for the structure;
 \item[{\bf A2.}] A fluid sub-problem, coupled with structure inertia (and the viscous part of the structure equation if $\epsilon \ne 0$).
 \end{enumerate}
 We will be numerically implementing this splitting by further splitting the fluid sub-problem into a dissipative and a
 non-dissipative part. Namely, we will be using the Lie splitting to separate the time-dependent Stokes
 problem from a pure advection problem. This way one can use non-dissipative solvers for non-dissipative problems, thereby 
 achieving higher accuracy. This is particularly relevant if one is interested in solving advection-dominated problems, such as, 
 for example, transport of nano-particles by blood. However, this additional splitting is not necessary to achieve stability of the scheme,
 as we show in Section~\ref{stability}, were only problems A1 and A2 above will be used to show 
 an energy estimate associated with unconditional stability of the scheme.
  However, for completeness, we show in the next section the implementation of our numerical scheme
 in which  the fluid sub-problem is further
 split into two, leading to the following Lie splitting of the coupled FSI problem:
 
 \begin{enumerate}
  \item [ ] {\bf A1.\phantom{(a)}} An elastodynamics sub-problem for the structure;
 \item[ ]{\bf A2(a).} A time-dependent Stokes problem for the fluid; 
 \item[ ] {\bf A2(b).} A fluid and ALE advection problem.
 \end{enumerate}

\subsubsection{Details of the operator-splitting scheme}\label{sec:scheme}
Before we define each of the sub-problems mentioned above, 
 we split the fluid stress $\widehat{\mathbf{\boldsymbol\sigma n}}$  into two parts, Part I and Part II:
$$
\widehat{\mathbf{\boldsymbol\sigma n}} = \underbrace{\widehat{\mathbf{\boldsymbol\sigma n}}+\beta\widehat{p\mathbf{n}}}_{(I)}
\underbrace{- \beta \widehat{p\mathbf{n}}}_{(II)},
$$
where $\beta$ is a number between $0$ and $1$, $0 \le \beta \le 1$.
Different values of $\beta$ effect the accuracy of the scheme \cite{Martina_paper1}, but not stability \cite{SunBorMar}.
For the problem discussed in the current manuscript, our experimental observations
indicate that $\beta = 1$ provides the kinematically-coupled $\beta$-scheme with
highest accuracy.

In addition to the splitting of the fluid stress, 
we also separate different physical effects in the structure problem. 
We split the viscoelastic effects from the purely elastic effects, and treat the structural viscoelasticity together with the fluid,
while treating the pure elastodynamics of the structure in a separate, hyperbolic, structural sub-problem.
Details of the splitting are as follows:

\vskip 0.1in
\noindent
\textbf{Problem A1:}  A pure elastodynamics problem is solved on $\hat\Omega^s$ with a boundary condition on $\hat\Gamma$
involving Part II of the normal fluid stress. Thus, the structure elastodynamics is driven by the initial velocity of the structure
(given by the fluid velocity at the fluid-structure interface from the previous time step), 
and by the fluid pressure loading $\beta  J  \widehat{(p {\bf n}^f)}$ acting on
the fluid-structure interface (obtained from the previous time step).
The problem reads:
Find $\boldsymbol{v}$, $\hat{\boldsymbol U}$ and $\hat{\boldsymbol V}$, with $p^{n}$  and $J^n$ obtained at the previous time step, such that
\begin{equation*}
 \left\{\begin{array}{l@{\ }} 
\displaystyle{\frac{\partial \boldsymbol{v}}{\partial t}} \bigg|_{\hat{\boldsymbol x}} = 0,   \quad \textrm{in} \; \Omega^f(t^n)\times(t^n, t^{n+1}), \\ \\
\displaystyle{\rho_s \frac{\partial \hat{\boldsymbol V}}{\partial t}} + \gamma \hat{\boldsymbol U}= \nabla \cdot \boldsymbol S( \hat{\boldsymbol U}),  \quad  
\displaystyle{ \frac{\partial \hat{\boldsymbol U}}{\partial t}} (z,t)  = \hat{\boldsymbol V}, \quad \textrm{in} \; \hat{\Omega}^s \times(t^n, t^{n+1}),\\ \\
 - J^n \beta (\widehat{p^n \boldsymbol n^f})|_{\hat\Gamma} +  \boldsymbol {Sn}^s|_{\hat{\Gamma}} = 0,\  {\rm on} \ \hat\Gamma\times(t^n, t^{n+1}),
 \end{array} \right.  
\end{equation*}
with boundary conditions: 
\begin{equation*}
\hat{\boldsymbol U}|_{z=0,L} = 0,
\end{equation*}
\begin{equation*}
 \hat{U}_z =  0, \quad \boldsymbol n^s_{ext} \cdot \boldsymbol {S n}^s_{ext} =  -P_{ext} \quad \textrm{on} \; \hat{\Gamma}^s_{ext} \times (t^n,t^{n+1}),
\end{equation*}
and initial conditions:
$$
\boldsymbol{v}(t^n)=\boldsymbol{v}^{n}, \quad \hat{\boldsymbol U}(t^n)=\hat{\boldsymbol U}^{n}, \quad \hat{\boldsymbol V}(t^n)=\hat{\boldsymbol V}^{n}.
$$
Then set $\boldsymbol{v}^{n+1/3}=\boldsymbol{v}(t^{n+1}), \; \hat{\boldsymbol U}^{n+1/3}=\hat{\boldsymbol U}(t^{n+1}), \; \hat{\hat{\boldsymbol V}}^{n+1/3}=\hat{\boldsymbol V}(t^{n+1}).$ 

A new ALE velocity $\boldsymbol{w}^{n+1}$ is calculated based on the current and previous locations of the fluid domain:
$\boldsymbol{w}^{n+1} = \partial \mathcal{A}_t/\partial t\ (\approx (\mathcal{A}^{n+1}_t-\mathcal{A}^n_t)/\Delta t$).
\vskip 0.1in
\noindent

\textbf{Problem A2(a).} A time-dependent Stokes problem is solved on the fixed fluid domain $\Omega^f(t^{n})$,
coupled with the viscous part of the structure problem defined on $\hat{\Omega}^s$.
(For higher accuracy, $\Omega^f(t^n)$ can be replaced with $\Omega^f(t^{n+1})$.)
The coupling is done via the kinematic coupling condition and a portion of
the dynamic coupling condition involving only Part I of the fluid stress. 
When $\epsilon = 0$ (the purely elastic case),
the problem on $\hat\Omega^s$ consists of only setting the structural velocity $\hat{\boldsymbol V}$ equal to the structural velocity from the 
previous time step, since $\partial\hat{\boldsymbol V}/\partial t = 0$ in this sub-problem.  The problem reads as follows:
 Find $\boldsymbol{v}, p, \hat{\boldsymbol V}$ and $\hat{\boldsymbol U}$, with $\hat{\boldsymbol{v}}(\hat{\boldsymbol{x}},t) = \boldsymbol{v}(\mathcal{A}_t(\hat{\boldsymbol{x}}),t)$,
such that for $t\in (t^n, t^{n+1})$, with $p^n$ and $J^n$ obtained in the previous time step, the following holds:
\begin{equation}\label{step1}
 \left\{\begin{array}{l@{\ }}
 \rho_f \displaystyle{\frac{\partial \boldsymbol{v}}{\partial t}}\bigg|_{\hat{\boldsymbol x}} =\nabla \cdot \boldsymbol{\sigma},   \quad \nabla \cdot \boldsymbol{v}=0, \quad \textrm{in} \; \Omega^f(t^{n}) \times(t^n, t^{n+1}), \\ \\
\rho_s \displaystyle{\frac{\partial \hat{\boldsymbol V}}{\partial t}} =  \epsilon \Delta \hat{\boldsymbol V}, \quad 
\displaystyle{ \frac{\partial \hat{\boldsymbol U}}{\partial t}} (\hat{z},t)  = 0, \quad \textrm{in} \; \hat{\Omega}^s \times(t^n, t^{n+1}), \\ \\ 
 J^{n+1} \ \widehat{\boldsymbol{\sigma n}^f}|_{\hat\Gamma} + J^n \beta (\widehat{p^n \boldsymbol n^f})|_{\hat\Gamma} + 
 \displaystyle{\epsilon \frac{\partial \hat{\boldsymbol V}}{\partial \boldsymbol n^s}\big|_{\hat{\Gamma}}}=0, \ {\textrm{on}}\ \hat\Gamma\times (t^n,t^{n+1}),\\ \\
\hat{\boldsymbol V}|_{\hat{\Gamma}} = \hat{\boldsymbol v}|_{\hat\Gamma}, \ {\textrm{on}}\ \hat\Gamma\times (t^n,t^{n+1}),
 \end{array} \right.  \nonumber
\end{equation}
with the following boundary conditions on $\Gamma^f_{\rm in}\cup\Gamma^f_{\rm out}\cup\Gamma^f_b$: 
\begin{equation*}
  \frac{\partial v_z}{\partial r}(z,0,t) = \quad v_r(z,0,t) = 0 \quad \textrm{on} \; \Gamma^f_b,
\end{equation*}
\begin{equation*}
  \boldsymbol{v}(0,R,t) = \boldsymbol{v}(L,R,t) = 0, 
\end{equation*}
\begin{equation*}
\boldsymbol{\sigma n}^f_{in} = -p_{in}(t)\boldsymbol{n}^f_{in}\  {\rm on}\ \Gamma^f_{\rm in}, \; \; \boldsymbol{\sigma n}^f_{out} = -p_{out}(t)\boldsymbol{n}^f_{out}  \ {\rm on}\ 
 \Gamma^f_{\rm out},
\end{equation*}
and initial conditions
$$
\boldsymbol{v}(t^n)=\boldsymbol{v}^{n+1/3}, \quad \hat{\boldsymbol U}(t^n)=\hat{\boldsymbol U}^{n+1/3}, \quad \hat{\boldsymbol V}(t^n)=\hat{\boldsymbol V}^{n+1/3}.
$$
Then set $\boldsymbol{v}^{n+2/3}=\boldsymbol{v}(t^{n+1}), \; \hat{\boldsymbol U}^{n+2/3}=\hat{\boldsymbol U}(t^{n+1}),\; \hat{\boldsymbol V}^{n+2/3}=\hat{\boldsymbol V}(t^{n+1}), \; p^{n+1}=p(t^{n+1}).$ 

\vskip 0.1in
\noindent
\textbf{Problem A2(b):} A fluid and ALE advection problem is solved on a fixed fluid domain 
 $\Omega^f(t^n)$ with the ALE velocity $\boldsymbol{w}^{n+1}$ calculated in Problem A1. 
 (For higher accuracy $\Omega(t^n)$ can be replaced by $\Omega(t^{n+1})$.) The problem reads:
Find $\boldsymbol{v}$, $\hat{\boldsymbol U}$ and $\hat{\boldsymbol V}$ with $\hat{\boldsymbol{v}}(\hat{\boldsymbol{x}},t) = \boldsymbol{v}(\mathcal{A}_t(\hat{\boldsymbol{x}}),t)$,
such that for $t\in (t^n, t^{n+1})$
\begin{equation*}
 \left\{\begin{array}{l@{\ }} 
\displaystyle{ \frac{\partial \boldsymbol{v}}{\partial t}}\bigg|_{\hat{\boldsymbol x}} + (\boldsymbol{v}^{n+2/3}-\boldsymbol{w}^{n+1}) \cdot \nabla \boldsymbol{v}= 0,   \quad \textrm{in} \; \Omega^f(t^n)\times (t^n, t^{n+1}),  \\ \\
\displaystyle{\frac{\partial \hat{\boldsymbol V}}{\partial t}} (\hat{z},t)  = 0, \quad
\displaystyle{\frac{\partial \hat{\boldsymbol U}}{\partial t}} (\hat{z},t)  = 0, \quad \textrm{in} \; \hat{\Omega}^s\times(t^n, t^{n+1}), \end{array} \right. 
\end{equation*}
with boundary conditions: 
$$\boldsymbol{v}=\boldsymbol{v}^{n+2/3} \  \; \textrm{on} \; \Gamma_{-}^{n+2/3}, \; \textrm{where}$$ 
$$\Gamma_{-}^{n+2/3} = \{\boldsymbol{x} \in \mathbb{R}^2 | \boldsymbol{x} \in \partial \Omega^f(t^n), (\boldsymbol{v}^{n+2/3}-\boldsymbol{w}^{n+1})\cdot \boldsymbol{n}^f <0 \},$$ 
and initial conditions
$$
\boldsymbol{v}(t^n)=\boldsymbol{v}^{n+2/3}, \quad \hat{\boldsymbol U}(t^n)=\hat{\boldsymbol U}^{n+2/3}, \quad \hat{\boldsymbol V}(t^n)=\hat{\boldsymbol V}^{n+2/3}.
$$
Then set $\boldsymbol{v}^{n+1}=\boldsymbol{v}(t^{n+1}), \; \hat{\boldsymbol U}^{n+1}=\hat{\boldsymbol U}(t^{n+1}), \; \hat{\boldsymbol V}^{n+1}=\hat{\boldsymbol V}(t^{n+1}).$

Do $t^n=t^{n+1}$ and return to Problem A1.

\noindent{\bf Remark.} The method proposed above works well even if the fluid and structure steps are performed in reverse order.

\subsection{Discretized subproblems}\label{algnum}
In this subsection we discretize each problem in space and time, and describe our solution strategy. 
To discretize the subproblems in time, we sub-divide the time interval $(0,T)$ into $N$ sub-intervals of width $\Delta t$.
For the space discretization, we use the finite element method. Thus, for $t^n = n \Delta t, 0\le n \le N,$ we define the finite
element spaces $V_h^f(t^n) \subset V^f(t^n), Q^f_h(t^n) \subset Q^f(t^n), V_h^{fsi}(t^n) \subset V^{fsi}(t^n)$ and $\hat{V}^s_h \subset \hat{V}^s$.
To write the weak formulation  the following notation for the corresponding bilinear forms will be used:
\begin{align}
 a_f(\boldsymbol v, \boldsymbol \varphi^f) &:= 2 \mu_f \int_{\Omega^f(t^n)} \boldsymbol D(\boldsymbol v) : \boldsymbol D(\boldsymbol \varphi^f) d \boldsymbol x,  \label{af}\\
 b_f(p_f, \boldsymbol \varphi^f) &:=  \int_{\Omega^f(t^n)} p_f \nabla \cdot \boldsymbol \varphi^f d \boldsymbol x, \label{bf} \\
  a_v(\hat{\boldsymbol V}, \hat{\boldsymbol \varphi}^v) &:= \epsilon \int_{\hat{\Omega}^s} \nabla \hat{\boldsymbol V} : \nabla \hat{\boldsymbol \varphi}^v d \boldsymbol x,  \label{av}\\
 a_e(\hat{\boldsymbol U}, \hat{\boldsymbol \varphi}^s) &:= 2 \mu_s \int_{\hat{\Omega}^s} \boldsymbol D(\hat{\boldsymbol U}) : \boldsymbol D(\hat{\boldsymbol \varphi}^s) d \boldsymbol x + \lambda_s \int_{\hat{\Omega}^s} (\nabla \cdot \hat{\boldsymbol U})(\nabla \cdot \hat{\boldsymbol \varphi}^s)  d \boldsymbol x    \nonumber \\
 & \quad +\gamma \int_{\hat{\Omega}^s} \hat{\boldsymbol U} \cdot \hat{\boldsymbol \varphi}^s  d \boldsymbol x  \label{ae}
 \end{align}
 
\noindent
\textbf{Problem A1:}
To discretize problem A1 in time we used the second order Newmark scheme. The weak formulation of the fully discrete problem is given as follows:
Find $\hat{\boldsymbol U}_h^{n+1/3} \in \hat{V}_h^s$ and $\hat{\boldsymbol V}_h^{n+1/3} \in \hat{V}_h^s$ such that for all $(\hat{\boldsymbol \varphi}^s_h,\boldsymbol \phi^s_h) \in \hat{V}_h^s \times \hat{V}_h^s$, with $p_h^n$ obtained at the previous time step:
\begin{equation*}
\begin{array}{rcl}
\displaystyle{\rho_{s} \int_{\hat{\Omega}^s} \frac{\hat{\boldsymbol V}_h^{n+1/3}-\hat{\boldsymbol V}_h^{n}}{\Delta t}\cdot \hat{\boldsymbol \varphi}_h^s d \boldsymbol x + a_e(\frac{\hat{\boldsymbol U}_h^{n}+\hat{\boldsymbol U}_h^{n+1/3}}{2}, \hat{\boldsymbol \varphi}_h^s)}&=&J^n \beta \displaystyle\int_{\hat{\Gamma}} \widehat{p_h^n \boldsymbol n^f} \cdot \hat{\boldsymbol \varphi}_h^s d \boldsymbol x,
\\
\displaystyle{ \rho_{s} \int_{\hat{\Omega}^s} (\frac{\hat{\boldsymbol V}_h^{n}+\hat{\boldsymbol V}_h^{n+1/3}}{2}- \frac{\hat{\boldsymbol U}_h^{n+1/3}-\hat{\boldsymbol U}_h^{n}}{\Delta t} ) \cdot \boldsymbol \phi_h^s d \boldsymbol x} &=& 0. 
\end{array}
\end{equation*}

\noindent
\textbf{Problem A2(a):}
We discretize problem A2(a) in time using the Backward Euler method. The weak formulation of the fully discrete problem is given as follows:
Find $(\boldsymbol{v}_h^{n+2/3},\hat{\boldsymbol V}_h^{n+2/3},p_{h}^{n+2/3}) \in V_h^{fsi}(t^n) \times Q^f_h(t^n)$ such that for all $(\boldsymbol \varphi_h^f,\hat{\boldsymbol \varphi}_h^v,\psi_h^f) \in V_h^{fsi}(t^n) \times Q_h^f(t^n)$:
\begin{equation*}
\rho_f \int_{\Omega^f(t^n)} \frac{\boldsymbol v^{n+2/3}_h-\boldsymbol v^{n+1/3}_h}{\Delta t} \cdot \boldsymbol \varphi^f_h d\boldsymbol x+
a_f(\boldsymbol v_h^{n+2/3}, \boldsymbol \varphi^f_h)- b_f(p^{n+2/3}_{h}, \boldsymbol \varphi^f_h)
+  b_f(\psi^f_h, \boldsymbol v_h^{n+2/3}) 
\end{equation*}
\begin{equation*}
 +\rho_s \int_{\hat{\Omega}^s} \frac{\hat{\boldsymbol V}_h^{n+2/3}-\hat{\boldsymbol V}^{n+1/3}_h}{\Delta t} \cdot \hat{\boldsymbol \varphi}^v_h d\boldsymbol x+  a_v(\hat{\boldsymbol V}_h^{n+2/3}, \hat{\boldsymbol \varphi}_h^v)
  = \int_0^R p_{in}(t^{n+1}) \varphi^f_{x,h}|_{x=0} dy
 \end{equation*}
 \begin{equation*}
 - \int_0^R p_{out}(t^{n+1}) \varphi^f_{x,h}|_{x=L} dy -J^n \beta \displaystyle\int_{\Gamma(t^n)} p_h^n \boldsymbol n^f \cdot \boldsymbol \varphi_h^f d \boldsymbol x.
 \end{equation*}
Then set $p_h^{n+1}=p_h^{n+2/3}$.
 
 \noindent
\textbf{Problem A2(b):}
We discretize problem A2(b) in time using again the Backward Euler method. To solve the resulting advection problem, we use a positivity-preserving ALE finite element scheme,
which preserves conservation of mass at the discrete level. Details of the scheme are presented in~\cite{boiarkine2011positivity}.

\section{Energy estimate associated with unconditional stability of the scheme for $\beta=0$}\label{stability}
We will show that the proposed operator splitting scheme satisfies an energy estimate which 
is associated with unconditional stablity of the scheme for all the parameters in the problem when $\beta = 0$. 
This will be done for the basic splitting algorithm, mentioned in Section~\ref{LieSplittingGeneral}, where the splitting consists of solving 
two sub-problems:
\begin{enumerate}
 \item[{\bf A1.}] An elastodynamics sub-problem for the structure;
 \item[{\bf A2.}] A fluid sub-problem, coupled with structure inertia (and the viscous part of the structure equation if $\epsilon \ne 0$).
 \end{enumerate}
We will use energy estimates to show that the energy of the discretized problem mimics the energy of
the continuous problem. More precisely, 
we will show that the total energy of the discretized problem plus 
viscous dissipation of the discretized problem, 
are bounded by the discrete energy of the initial data and the work done by the inlet and outlet dynamic pressure data.
In contrast with similar results appearing in literature \cite{causin2005added,SunBorMar,badia,Fernandez3,burman2009stabilization}, which consider simplified models without fluid advection, 
and linearized fluid-structure coupling calculated at a fixed fluid domain boundary, 
in this manuscript we derive the corresponding energy estimate for a full, nonlinear FSI problem, 
that includes fluid advection, and the coupling is achieved at the moving fluid-structure interface.

To simplify analysis, the following assumptions that do not influence stability of the scheme, will be considered:
\begin{description}
\item[1.] Only radial displacement of the fluid-structure interface is allowed, i.e., $\hat{\boldsymbol U}|_{r=R}\cdot \hat{\boldsymbol e}_z = 0$
 at the fluid-structure interface. The FSI problem with this boundary condition is well-defined. This assumption does not affect  stability of the scheme related to the added mass effect.
In fact, the same assumption was considered in the original work on added mass effect by Causin et al.~\cite{causin2005added}.
In the present manuscript, this simplifies the form of the energy estimates in the proof.
\item[2.] The problem is driven by the dynamic inlet and outlet pressure data, and the flow enters and leaves
the fluid domain parallel to the horizontal axis: 
$$
p+\frac{\rho_f}{2} |\boldsymbol v|^2 = p_{in/out}(t), \  v_r = 0, \ {\rm on}\ \Gamma_{in/out}.
$$
\end{description}

We consider $\beta = 0$ here because, in this case, it is easier to prove the related energy estimates.
Our numerical results presented in \cite{Martina_paper1} indicate that only accuracy, not stability, is affected by changing
$\beta$ between $0$ and $1$.
Using energy estimates to prove unconditional stability of the scheme with $\beta \ne 0$ would be significantly more difficult.
We mention a related work, however, in which unconditional stability of the kinematically-coupled $\beta$-scheme
for $\beta\in[0,1]$ was proved for a simplified, linearized FSI problem with a thin structure, using different techniques from those
presented here \cite{SunBorMar}.
In the present paper, for the first time, we derive an energy estimate associated with unconditional stability   
of the full, nonlinear FSI problem, defined on a moving domain, with nonlinear fluid-structure coupling.
The same proof applies to problems where the thick structure is replaced by a thin structure.
In that case, the kinematically coupled $\beta$-scheme is a fully partitioned scheme \cite{guidoboni2009stable,Martina_paper1}.

We begin by first deriving an energy equality of the continuous, coupled FSI problem.

\subsection{The energy of the continuous coupled problem}
To formally derive an energy equality of the coupled FSI problem we multiply the structure equations by the structure velocity,
the balance of momentum in the fluid equations by the fluid velocity, and integrate by parts over the respective domains
using the incompressibility condition.
The dynamic and kinematic coupling conditions are then used to couple the fluid and structure sub-problems.
The resulting equation represents the total energy of the problem.

The following identities are used in this calculation:
\begin{eqnarray}
\int_{\Omega^f(t)} \frac{\partial \boldsymbol v}{\partial t} \boldsymbol v d\boldsymbol x &=& \frac{1}{2} \frac{d}{dt}\int_{\Omega^f(t)} |\boldsymbol v|^2 d\boldsymbol x - 
\frac{1}{2} \int_{\Gamma(t)} |\boldsymbol v|^2 \boldsymbol v \cdot \boldsymbol n^f dS, \label{ident1} 
\\
\int_{\Omega^f(t)} (\boldsymbol v \cdot \nabla) \boldsymbol v \cdot \boldsymbol v d\boldsymbol x &=& \frac{1}{2} \int_{\partial \Omega^f(t)} |\boldsymbol v|^2 \boldsymbol v \cdot \boldsymbol n^f dS, \label{ident2}
\end{eqnarray}
The first  is just the transport theorem. 
The second one is obtained using integration by parts.

The boundary integral associated with the convective term \eqref{ident2} is simplified as follows.
The portion corresponding to $\Gamma_b$ is zero due to the 
symmetry boundary condition, which implies ${\bf v}\cdot{\bf n}=0$ on $\Gamma_b$.
The portion corresponding to $\Gamma(t)$ is canceled with the same term appearing in the transport formula \eqref{ident1}. Finally,
the boundary terms on $\Gamma_{in/out}$ are absorbed by the dynamic pressure boundary conditions. 
The fluid sub-problem implies:
\begin{eqnarray*}
&\displaystyle{ \frac{1}{2}\frac{d}{dt} \bigg\{  \rho_f ||\boldsymbol v||^2_{L^2(\Omega^f(t))} \bigg\} + 2 \mu_f ||\boldsymbol D(\boldsymbol v)||^2_{L^2(\Omega(t))}}\\
&\displaystyle{ - \int_0^R p_{in}(t) v_z|_{z=0} dr +   \int_0^R p_{out}(t) v_z|_{z=L} dr = \int_{\Gamma(t)} \boldsymbol\sigma \boldsymbol {n}^f \cdot \boldsymbol v\ dS}
\end{eqnarray*}
The integral on the right-hand side can be written in Lagrangian coordinates as
\begin{equation}
\int_{\Gamma(t)} \boldsymbol\sigma \boldsymbol{n}^f \cdot \boldsymbol v\ dS
=\int_{\hat{\Gamma}}\widehat{\boldsymbol\sigma \boldsymbol{n}^f} \cdot \hat{\boldsymbol v} \ J \ d\hat{z}
\end{equation}
where $J$ is the Jacobian of the transformation from the Eulerian to Lagrangian coordinates.
We use the kinematic and dynamic  lateral boundary conditions \eqref{kinematic}-\eqref{dynamic} to obtain
\begin{equation}\label{weak_coupling}
\int_{\hat{\Gamma}} \widehat{\boldsymbol\sigma \boldsymbol {n}^f} \cdot \hat{\boldsymbol v}  \ J \ d\hat{z}
= -  \int_{\hat{\Gamma}} \left[\boldsymbol S \boldsymbol n^s  \cdot \frac{\partial \hat{\boldsymbol U}}{\partial t}-\epsilon \frac{\partial }{\partial \boldsymbol n^s}\bigg(\frac{\partial \hat{\boldsymbol U}}{\partial t} \bigg)  \cdot \frac{\partial \hat{\boldsymbol U}}{\partial t} \right]\bigg|_{\hat{\Gamma}}  d\hat{z}.
\end{equation}
After adding the energy equalities associated with the fluid problem and the thick structure problem,
and after using the coupling expressed in \eqref{weak_coupling}, one obtains the following energy equality of the coupled FSI problem:
 \begin{eqnarray}
&\displaystyle{ \frac{d}{dt} \Bigg\{ \underbrace{ \frac{\rho_f}{2} ||\boldsymbol v||^2_{L^2(\Omega^f(t))} 
+  \frac{\rho_s}{2} \bigg|\bigg|\frac{\partial \hat{\boldsymbol U}}{\partial t} \bigg|\bigg|^2_{L^2(\hat{\Omega}^s)}}_{Fluid\ and \ Structure \ Kinetic \ Energy}}
 \nonumber\\
&  \underbrace{+ \frac{\gamma}{2}  ||\hat{\boldsymbol U}||^2_{L^2(\hat{\Omega}^s)}+\mu_s  ||\boldsymbol D(\hat{\boldsymbol U})||^2_{L^2(\hat{\Omega}^s)}+\frac{\lambda_s}{2} ||\nabla \cdot \hat{\boldsymbol  U}||^2_{L^2(\hat{\Omega}^s)} }_{Structure \ Elastic \ Energy\ }
\Bigg\}
 \nonumber
\\
 & +  \underbrace{2 \mu_f ||\boldsymbol D(\boldsymbol v)||^2_{L^2(\Omega^f(t))}}_{Fluid \ Viscous \ Energy}+  \underbrace{\epsilon \bigg|\bigg|\nabla \frac{\partial \hat{\boldsymbol U}}{\partial t}\bigg|\bigg|^2_{L^2(\hat{\Omega}^s)}}_{Structure \ Viscous \ Energy}
 \nonumber
\\
& = 
  \displaystyle{\int_0^R p_{in}(t) v_z|_{z=0} dr -   \int_0^R p_{out}(t) v_z|_{z=L} dr -  \int_{\hat{\Gamma}^s_{ext}} P_{ext} \frac{\partial \hat{U}_r}{\partial t} dS}.
 \label{FSIEnergy}
\end{eqnarray}

To obtain energy estimates for the proposed operator splitting scheme, we first write the main steps of the splitting scheme
in weak form. For this purpose we start by writing the weak form of the problem written in ALE formulation, with dynamics pressure
inlet and outlet data, and
then split the weak ALE form, following the operator splitting approach presented in the previous section.
 
\subsection{The weak form of the continuous coupled problem in ALE form}
We consider the ALE form of the fluid equations \eqref{ALEsys1}-\eqref{ALEsys2}, coupled with the initial and
boundary conditions \eqref{inlet}-\eqref{dynamic}, where the inlet and outlet conditions for the fluid problem are
the {\sl dynamic pressure data}.
 As we shall see later, it is convenient to write the fluid and ALE advection term in symmetric form, 
 giving rise to the following weak formulation:
for $t \in (0, T)$, find $(\boldsymbol v, \hat{\boldsymbol U}, p) \in V^{f}(t)\times \hat{V}^s \times Q(t)$ such that 
$\boldsymbol v = {\hat{\boldsymbol U}}_t \circ \mathcal A_t^{-1}$ on $\Gamma(t)$, and the following holds:
\begin{equation*}
\rho_f \int_{\Omega^f(t)}  \frac{\partial \boldsymbol{v}}{\partial t}\bigg|_{\hat{\boldsymbol x}}\cdot \boldsymbol \varphi d\boldsymbol x+ 
\frac{\rho_f}{2}\int_{\Omega^f(t)}\Big ( ((\boldsymbol v - {\bf w})\cdot \nabla) \boldsymbol v \cdot \boldsymbol \varphi
-((\boldsymbol v - {\bf w})\cdot \nabla) \boldsymbol \varphi \cdot \boldsymbol v \Big )d \boldsymbol x
\end{equation*}  
\begin{equation*}
+\frac{\rho_f}{2}\int_{\Omega^f(t)}(\nabla\cdot{\bf w}){\bf v}\cdot\boldsymbol\varphi d \boldsymbol x
+2 \mu_f \int_{\Omega^f(t)} \boldsymbol D(\boldsymbol v) : \boldsymbol D(\boldsymbol\varphi)d \boldsymbol x 
\end{equation*}
\begin{equation*}
- \int_{\Omega^f(t)} p \nabla \cdot \boldsymbol \varphi d \boldsymbol x +\rho_s \int_{\hat{\Omega}^s} \frac{\partial^2 \hat{\boldsymbol U}}{\partial t^2} \cdot \hat{ \boldsymbol \psi} d\boldsymbol x+2 \mu_s \int_{\hat{\Omega}^s} \boldsymbol D(\hat{\boldsymbol U}) : \boldsymbol D(\hat{ \boldsymbol \psi}) d \boldsymbol x 
\end{equation*}
\begin{equation*}
+\lambda_s \int_{\hat{\Omega}^s} (\nabla \cdot \hat{\boldsymbol U}) (\nabla \cdot \hat{ \boldsymbol \psi}) d \boldsymbol x +\gamma \int_{\hat{\Omega}^s} \hat{\boldsymbol U} \cdot \hat{ \boldsymbol \psi} d\boldsymbol x+ \epsilon \int_{\hat{\Omega}^s} \nabla \frac{\partial \hat{\boldsymbol U}}{\partial t} : \nabla \hat{ \boldsymbol \psi} d\boldsymbol x
\end{equation*}
\begin{equation}
 + \int_{\hat{\Gamma}^s_{ext}} P_{ext} (\hat{ \boldsymbol \psi} \cdot \boldsymbol n^s_{ext}) dS= \int_0^R p_{in}(t) \varphi_z|_{z=0} dr-  \int_0^R p_{out}(t) \varphi_z|_{z=L} dr,\label{weak_symmetrized}
\end{equation}
and 
\begin{equation*}
\int_{\Omega^f(t)} q \nabla \cdot \boldsymbol v d \boldsymbol x = 0, 
\end{equation*}
for all $(\boldsymbol \varphi, \hat{ \boldsymbol \psi}, q) \in V^{fsi}(t) \times Q(t)$,
where $V^{fsi}(t)$ and $Q(t)$ are defined by \eqref{Vf(t)}-\eqref{Vfsi}.

This weak formulation is consistent with the problem. Indeed, integration by parts of one half of the convective term gives:
$$
\frac{1}{2}\int_{\Omega^f(t)} ((\boldsymbol v - {\bf w})\cdot \nabla) \boldsymbol v \cdot \boldsymbol \varphi=
-
\frac{1}{2}\int_{\Omega^f(t)} ((\boldsymbol v - {\bf w})\cdot \nabla) \boldsymbol \varphi \cdot \boldsymbol v
+\frac{1}{2} \int_{\Omega^f(t)} (\nabla\cdot{\bf w}){\bf v}\cdot \boldsymbol \varphi
$$
$$
+\frac{1}{2}\int_{\partial\Omega^f(t)}((\boldsymbol v - {\bf w})\cdot {\bf n}) \boldsymbol v \cdot \boldsymbol \varphi.
$$
Here, we have used the fact that $\nabla\cdot{\bf v}=0$. The last term in this expression, i.e., the boundary integral, 
can be evaluated as follows:
on $\Gamma(t)$ we recall 
that ${\bf v}={\bf w}$ and so that part is zero; on $\Gamma_b$ we have ${\bf v}\cdot{\bf n}={\bf w}\cdot{\bf n}=0$, and so this contribution
is zero as well; finally, on $\Gamma_{in/out}$ we have ${\bf w} = 0$ and the remaining quadratic velocity term is exactly 
the quadratic velocity contribution in dynamic pressure. Therefore, the weak form \eqref{weak_symmetrized} is consistent
with problem \eqref{ALEsys1}-\eqref{ALEsys2} and
 \eqref{inlet}-\eqref{dynamic} in ALE form, with dynamic inlet and outlet pressure data. 

\subsection{The time discretization via operator splitting in weak form}
We perform the time discretization via operator splitting, described in Section~\ref{sec:scheme}, 
and write each of the split sub-problems in weak form.
To simplify the energy estimate proof,  we keep the entire fluid sub-problem together, without splitting the advection problem from the time-dependent Stokes problem. 
The main features of the scheme, which are related to how the fluid and structure problems are separated, are left intact in
the stability analysis.
Therefore, we split the coupled FSI problem into:
\begin{enumerate}
 \item[{\bf A1.}] An elastodynamics sub-problem for the structure;
 \item[{\bf A2.}] A fluid sub-problem, coupled with structure inertia (and the viscous part of the structure equation if $\epsilon \ne 0$).
 \end{enumerate}

We discretize the problem in time as described in Section~\ref{algnum}. Since our stability analysis does not depend on the spatial discretization, 
we leave the spatial operators in continuous form. 
To write the weak forms of the structure and fluid sub-problems we use the bilinear forms defined in~\eqref{af}-\eqref{ae}.
A weak formulation of our semi-discrete operator splitting numerical scheme is given as follows:
\vskip 0.1in
 \textbf{Problem A1. (The structure sub-problem)} Find $\hat{\boldsymbol U}^{n+1/2} \in \hat{V}^s$ and $\hat{\boldsymbol V}^{n+1/2} \in \hat{V}^s$ such that for all $(\hat{\boldsymbol \varphi}^s,\boldsymbol \phi^s) \in \hat{V}^s \times \hat{V}^s$ and $t\in(t^n,t^{n+1}):$ 
\begin{equation}\label{S2discrete}
\begin{array}{rcl}
\displaystyle{\rho_{s} \int_{\hat{\Omega}^s} \frac{\hat{\boldsymbol V}^{n+1/2}-\hat{\boldsymbol V}^{n}}{\Delta t}\cdot \hat{\boldsymbol \varphi}^s d \boldsymbol x + a_e(\frac{\hat{\boldsymbol U}^{n}+\hat{\boldsymbol U}^{n+1/2}}{2}, \hat{\boldsymbol \varphi}^s)}&=&0,
\\
\displaystyle{ \rho_{s} \int_{\hat{\Omega}^s} (\frac{\hat{\boldsymbol V}^{n}+\hat{\boldsymbol V}^{n+1/2}}{2}- \frac{\hat{\boldsymbol U}^{n+1/2}-\hat{\boldsymbol U}^{n}}{\Delta t} ) \cdot \boldsymbol \phi^s d \boldsymbol x} &=& 0. 
\end{array}
\end{equation}
 In this step $\partial {\boldsymbol v}/{\partial t} = 0,$ so $\boldsymbol v^{n+1/2}=\boldsymbol v^{n}.$
 
 \vskip 0.1in
 \textbf{Problem A2. (The fluid sub-problem)} Find $(\boldsymbol{v}^{n+1},\hat{\boldsymbol V}^{n+1},p^{n+1}) \in V^{fsi}(t) \times Q^f(t)$ such that for all $(\boldsymbol \varphi^f,\hat{\boldsymbol \varphi}^v,\psi^f) \in V^{fsi}(t) \times Q^f(t)$ and  $t\in (t^n, t^{n+1})$:
 
\begin{equation*}
\rho_f \int_{\Omega^f(t^n)} \frac{\boldsymbol v^{n+1}-\boldsymbol v^{n+1/2}}{\Delta t} \cdot \boldsymbol \varphi^f d\boldsymbol x+
\frac {\rho_f}{2} \int_{\Omega^f(t^n)}(\nabla\cdot{\boldsymbol w}^{n+\frac 1 2}){\boldsymbol v}^{n+1}\boldsymbol \varphi^f
\end{equation*}
\begin{equation*}
+\frac{\rho_f}{2}\int_{\Omega^f(t^n)}\left (({\boldsymbol v}^n-{\boldsymbol w}^{n+1/2})\cdot\nabla){\boldsymbol v}^{n+1}\cdot\varphi^f
-(({\boldsymbol v}^n-{\boldsymbol w}^{n+1/2})\cdot\nabla)\varphi^f\cdot{\boldsymbol v}_n^{n+1}\right )
\end{equation*}
\begin{equation*}
+ a_f(\boldsymbol v^{n+1}, \boldsymbol \varphi^f)- b_f(p^{n+1}, \boldsymbol \varphi^f)
 +\rho_s \int_{\hat{\Omega}^s} \frac{\hat{\boldsymbol V}^{n+1}-\hat{\boldsymbol V}^{n+1/2}}{\Delta t} \cdot \hat{\boldsymbol \varphi}^v d\boldsymbol x+  a_v(\hat{\boldsymbol V}^{n+1}, \hat{\boldsymbol \varphi}^v)
 \end{equation*}
 \begin{equation}
 = \int_0^R p_{in}(t^{n+1}) \varphi^f_{x}|_{x=0} dy - \int_0^R p_{out}(t^{n+1}) \varphi^f_{x}|_{x=L} dy, \label{S1discrete}
 \end{equation}
 \begin{equation*}
  b_f(\psi^f, \boldsymbol v^{n+1}) = 0. 
 \end{equation*}
 In this step $\partial {\hat{\boldsymbol U}}/{\partial t} = 0,$ so $\hat{\boldsymbol U}^{n+1}=\hat{\boldsymbol U}^{n+1/2}.$ 
 
 Notice that in Problem A2 we have taken the
 ALE velocity $\boldsymbol{w}^{n+1/2}$  from the just-calculated Problem A1,
 and not from the previous time step. Since in Problem A2 the fluid domain does not change,
 $\boldsymbol{w}^{n+1/2} = \boldsymbol{w}^{n+1}$.
 This is important in proving energy estimates that are associated with the stability of the splitting scheme for the full,
 nonlinear FSI problem defined on the moving domain. Namely, as we shall see below,
 by using this ALE velocity we will be able to 
approximate the total discrete energy at $t^{n+1}$, which includes the kinetic energy
 due to the motion of the fluid domain at $t^{n+1}$, described by $\boldsymbol{w}^{n+1}$.

\subsection{Energy estimate associated with unconditional stability of the splitting scheme}

Let $\mathcal{E}_f^n$ denote the discrete energy associated with~\eqref{S1discrete}, and let $\mathcal{E}_s^n$ denote the discrete energy 
associated with~\eqref{S2discrete} at time level $n \Delta t$:
\begin{eqnarray}
 \mathcal{E}_f^n & :=& \frac{\rho_f}{2} ||\boldsymbol v^n||^2_{L^2(\Omega^f(t^n))}, \\
 \mathcal{E}_s^n &:=& \frac{\rho_{s}}{2} ||\hat{\boldsymbol V}^n||^2_{L^2(\hat{\Omega}^s)} + \mu_s ||D(\hat{\boldsymbol U}^n)||^2_{L^2(\hat{\Omega^s})} + \frac{\lambda_s}{2}||\nabla \cdot \hat{\boldsymbol U}^n||^2_{L^2(\hat{\Omega}^s)} \nonumber\\
 & & \quad + \frac{\gamma}{2}||\hat{\boldsymbol U}^n||^2_{L^2(\hat{\Omega}^s)}.
\end{eqnarray}
The following energy estimate holds for the full nonlinear FSI problem,
satisfying assumptions {\bf 1} and {\bf 2} above:

\begin{theorem}{\bf (Energy estimate of the operator splitting scheme)}
Let $\{(\boldsymbol v^n,\hat{\boldsymbol V}^n, \hat{\boldsymbol U}^n \}_{0 \le n \le N}$ be a solution of~\eqref{S2discrete}-\eqref{S1discrete}.
Then, at any time-level $\tilde{n} \Delta t$, where $0\le \tilde{n} \le N$, the following energy estimate holds:
\begin{equation*}
 \mathcal{E}_f^{\tilde{n}}+\mathcal{E}_s^\tn 
+  \frac{\rho_f}{2} \sum_{n=0}^{\tn-1}  ||\boldsymbol v^{n+1}-\boldsymbol v^n||^2_{L^2(\Omega^f(t^n))} 
+ \frac{\rho_s}{2}  \sum_{n=0}^{\tn-1} ||\hat{\boldsymbol V}^{n+1}-\hat{\boldsymbol V}^{n+1/2}||^2_{L^2(\Omega^s)} 
 \end{equation*}
\begin{equation}\label{energy_estimate}
+ \mu_f \Delta t \sum_{n=0}^{\tn-1} ||\boldsymbol D(\boldsymbol v^{n+1})||^2_{L^2(\Omega^f(t^n))}
+ \epsilon \Delta t \sum_{n=0}^{\tn-1} ||\nabla \hat{\boldsymbol V}^{n+1}||^2_{L^2(\hat{\Omega}^s)}
\end{equation}
\begin{equation*}
  \leq  \mathcal{E}_f^0 + \mathcal{E}_s^0 +\frac{C \Delta t}{2 \mu_f} \sum_{n=0}^{\tn-1}||p_{in}(t^n)||^2_{L^2(0,R)}+\frac{C \Delta t}{2 \mu_f} \sum_{n=0}^{\tn-1}||p_{out}(t^n)||^2_{L^2(0,R)}.
\end{equation*}
 \end{theorem}
 The first line in the energy estimate corresponds to the kinetic energy of the fluid and the total energy of the structure, while
 the second line describes viscous dissipation in the fluid and the structure. They are estimated by
 the initial kinetic energy of the fluid, by the total initial energy of the structure, and by the work done
 by the inlet and outlet pressure data.
 \begin{proof}
To prove the energy estimate, we test the structure  problem \eqref{S2discrete} with 
$$
(\hat{\boldsymbol \varphi}^s, \hat{\boldsymbol \phi}^s)=\bigg(\frac{\hat{\boldsymbol U}^{n+1/2}-\hat{\boldsymbol U}^{n}}{\Delta t}, \frac{\hat{\boldsymbol V}^{n+1/2}-\hat{\boldsymbol V}^{n}}{\Delta t}\bigg),
$$
and problem ~\eqref{S1discrete} with 
$$
(\boldsymbol \varphi^f, \hat{\boldsymbol \varphi}^v, \psi^f) = (\boldsymbol v^{n+1}, \hat{\boldsymbol V}^{n+1}, p^{n+1}).
$$ 

Since we have  assumed that the fluid-structure interface deforms only in the radial direction, i.e. $\hat{\boldsymbol U}|_{r=R}=\hat{U}|_{r=R}{\hat{\boldsymbol e}}_r$,
we can explicitly calculate the ALE mapping at every step and, more importantly, calculate the associated ALE velocity ${\boldsymbol w}$.
For this purpose, we denote by $\eta$ the radial displacement of the fluid-structure interface,
namely,
$$
\hat{\eta} := \hat{\boldsymbol U}|_{\hat{r} = R} {\hat {\boldsymbol e}}_r, \ {\rm and}\  \hat{\eta}^n := \hat{\boldsymbol U}(\hat{z},t^n)|_{\hat{r} = R} {\hat{\boldsymbol e}}_r.
$$
We consider the following simple ALE mapping:
$$
A_{t^n}: \hat\Omega^f \to \Omega^f(t^n),\quad A_{t^n}(\hat{z},\hat{r}):=\left(\hat{z},\frac{R+\hat{\eta}^n}{R}\hat{r}\right)^\tau.
$$
We will  also need the explicit form of the ALE mapping from the computational 
domain $\Omega^f(t^n)$ to $\Omega^f(t^{n+1})$, which is given by
$$
A_{t^{n+1}} \circ A_{t^n}^{-1}: \Omega^f(t^n) \to \Omega^f(t^{n+1}), \quad A_{t^{n+1}} \circ A_{t^n}^{-1}(z,r) = \left(z, \frac{R + \hat{\eta}^{n+1}}{R+\hat{\eta}^n} r\right)^\tau.
$$
The corresponding Jacobian and the ALE velocity are given by
\begin{equation}\label{jacobian}
J_n^{n+1} := \frac{R + \hat{\eta}^{n+1}}{R+\hat{\eta}^n}, \quad 
{\boldsymbol w}^{n+1}
=\frac{1}{\Delta t}\frac{\hat\eta^{n+1}-\hat\eta^n}{R+\hat\eta^n} \ {r}\hat{{\boldsymbol e}}_r.
\end{equation}
In Problem A1 where we just calculated the updated location of the structure,
$\hat\eta^n$ determines the ``reference domain'', and $\hat\eta^{n+1/2}$, which is the same as $\hat\eta^{n+1}$, 
determines the location of the new domain. Therefore, 
the time-derivative of the interface displacement is approximated by $(\hat\eta^{n+1}-\hat\eta^n)/\Delta t$, which enters 
the expression for the ALE velocity ${\boldsymbol w}^{n+1}$. Again, notice that ${\boldsymbol w}^{n+1/2} = {\boldsymbol w}^{n+1}$.

We begin by considering the fluid sub-problem and the advection terms involving the fluid and ALE advection. 
After replacing the test functions by the fluid velocity at time $n+1$, we first notice that the symmetrized advection terms cancel out.
What is left are the terms 
$$
{\rho_f} \int_{\Omega^f(t^n)}  \frac{{\boldsymbol v}^{n+1} - {\boldsymbol v}^n}{\Delta t} {\boldsymbol v}^{n+1} 
+\frac{\rho_f}{2} \int_{\Omega^f(t^n)}(\nabla\cdot{\boldsymbol w}^{n+\frac 1 2})|{\boldsymbol v}^{n+1}|^2.
$$
To deal with the term on the left, we use the following identity:
\begin{equation}\label{famous_formula}
a(a-b) =\frac 1 2( a^2-b^2+(a-b)^2).
\end{equation}
To deal with the term on the right hand-side, we use the expression for ${\boldsymbol w}^{n+1/2}$, given by \eqref{jacobian},
where we recall that ${\boldsymbol w}^{n+1/2} = {\boldsymbol w}^{n+1}$, since in the second step, the structure location 
does not change. 
We obtain
$$
\frac {\rho_f}{2}\frac{1}{\Delta t}\int_{\Omega^f(t^n)}\left(1+\frac{\hat\eta^{n+1}-\hat\eta^n}{R+\hat\eta^n} \right)|{\boldsymbol v}^{n+1}|^2
+|{\boldsymbol v}^{n+1}-{\boldsymbol v}^{n}|^2-|{\boldsymbol v}^{n}|^2
$$
$$
=\frac{\rho_f}{2}\frac{1}{\Delta t}\int_{\Omega^f(t^n)}\Big (\frac{R+\hat\eta^{n+1}}{R+\hat\eta^n}|{\boldsymbol v}^{n+1}|^2
+|{\boldsymbol v}^{n+1}-{\boldsymbol v}^{n}|^2-|{\boldsymbol v}^{n}|^2\Big).
$$
Now notice that $({R+\hat\eta^{n+1}})/({R+\hat\eta^n})$ is exactly the Jacobian of the ALE mapping from $\Omega^f(t^n)$ to $\Omega^f(t^{n+1})$, 
see \eqref{jacobian}, and so we can convert that integral into an integral over $\Omega^f(t^{n+1})$ to recover
the kinetic energy of the fluid at the next time step:
$$
\frac{\rho_f}{2}\frac{1}{\Delta t}\int_{\Omega^f(t^n)}\frac{R+\hat\eta^{n+1}}{R+\hat\eta^n}|{\boldsymbol v}^{n+1}|^2
=\frac{\rho_f}{2}\frac{1}{\Delta t}\int_{\Omega^f(t^{n+1})}|{\boldsymbol v}^{n+1}|^2.
$$
This calculation effectively shows that the kinetic energy of the fluid at the next time step accounts for the 
kinetic energy of the fluid at the previous time step, plus the kinetic energy due to the motion of the fluid domain.
Therefore, even at the discrete level, we see that the energy of the discretized problem mimics well the energy of
the continuous problem.

Notice that this calculation also shows that the ALE mapping and its Jacobian satisfy the geometric conservation law,
see \cite{Farhat} for more details. A similar result was shown in \cite{LukacovaCMAME} for the ALE mapping which is
the harmonic extension of the boundary to the entire domain.

To deal with the structure sub-problem, we do not have the same problem associated with moving domains,
since the structure problem is defined
in Lagrangian coordinates, namely, on a fixed domain $\hat\Omega^s$. We use formula \eqref{famous_formula}
to calculate the kinetic energy of the structure at the next time step, in terms of the kinetic energy of the structure
at the previous time step, plus a quadratic term  
$||\hat{\boldsymbol V}^{n+1}-\hat{\boldsymbol V}^{n+1/2}||^2_{L^2(\hat{\Omega}^s)}$ that accounts for the kinetic energy
due to the difference in the velocities of the structure between the two time steps. 

Finally, we add the corresponding energy equalities for the fluid and structure sub-problems. 
We obtain:
\begin{equation*}
\mathcal{E}_f^{n+1}+\mathcal{E}_s^{n+1}+ \frac{\rho_f}{2}
||\boldsymbol v^{n+1}-\boldsymbol v^n||^2_{L^2(\Omega^f(t^n))} + \frac{\rho_s}{2}
||\hat{\boldsymbol V}^{n+1}-\hat{\boldsymbol V}^{n+1/2}||^2_{L^2(\hat{\Omega}^s)} 
 \end{equation*}
\begin{equation*}
 +2 \mu_f \Delta t ||D(\boldsymbol v^{n+1})||^2_{L^2(\Omega^f(t^n))} 
 + \epsilon \Delta t \sum_{n=0}^{\tn-1} ||\nabla \hat{\boldsymbol V}^{n+1}||^2_{L^2(\hat{\Omega}^s)}
\end{equation*}
\begin{equation*}
 = \mathcal{E}_f^n + \mathcal{E}_s^n +  \Delta t  \int_0^R p_{in}(t^{n+1}) \boldsymbol v^{n+1}|_{x=0} dy 
 - \Delta t \int_0^R p_{out}(t^{n+1}) \boldsymbol v^{n+1}|_{x=L} dy.
\end{equation*}
To bound the right-hand side of this equality, we use the Cauchy-Schwartz and Young's inequalities, to obtain:
\begin{equation*}
\Delta t  \int_0^R p_{in}(t^{n+1}) \boldsymbol v^{n+1}|_{x=0} dy  - \Delta t \int_0^R p_{out}(t^{n+1}) \boldsymbol v^{n+1}|_{x=L} dy 
\end{equation*}
\begin{equation*}
\le \frac{\Delta t}{2 \epsilon_1 } ||p_{in}(t^n)||^2_{L^2(0,R)}+\frac{\Delta t}{2 \epsilon_1 } ||p_{out}(t^n)||^2_{L^2(0,R)}+\epsilon_1 \Delta t
|| {\boldsymbol v}^{n+1}||^2_{L^2(0,R)}.
\end{equation*}
From the trace and Korn inequalities, we have
\begin{equation*}
\Delta t  \int_0^R p_{in}(t^{n+1}) \boldsymbol v^{n+1}|_{x=0} dy  - \Delta t \int_0^R p_{out}(t^{n+1})\boldsymbol v^{n+1}|_{x=L} dy 
\end{equation*}
\begin{equation*}
\le \frac{\Delta t}{2 \epsilon_1 } ||p_{in}(t^n)||^2_{L^2(0,R)}+\frac{\Delta t}{2 \epsilon_1 } ||p_{out}(t^n)||^2_{L^2(0,R)}+\epsilon_1 C \Delta t||\boldsymbol D( \boldsymbol v^{n+1})||^2_{L^2(\Omega^f(t^n))},
\end{equation*}
where $C$ is the constant from the trace and Korn inequalities. 
In general, Korn's constant depends on the domain. It was shown, however, that for domains associated 
with fluid-structure interaction problems of the type studied in this manuscript, the Korn's constant is independent
of the sequence of approximating domains \cite{BorSun,BorSunMulti}.

By setting $\epsilon_1 = \displaystyle\frac{\mu_f}{C}$, the last term can be combined with the term 
on the left hand-side, associated with fluid diffusion. 
Therefore, so far we have shown that the following inequality holds:
\begin{equation*}
\mathcal{E}_f^{n+1}+\mathcal{E}_s^{n+1}+ \frac{\rho_f}{2}
||\boldsymbol v^{n+1}-\boldsymbol v^n||^2_{L^2(\Omega^f(t^n))} + \frac{\rho_s}{2}
||\hat{\boldsymbol V}^{n+1}-\hat{\boldsymbol V}^{n+1/2}||^2_{L^2(\hat{\Omega}^s)} 
 \end{equation*}
\begin{equation*}
 +\mu_f \Delta t ||D(\boldsymbol v^{n+1})||^2_{L^2(\Omega^f(t^n))}
 + \epsilon \Delta t \sum_{n=0}^{\tn-1} ||\nabla \hat{\boldsymbol V}^{n+1}||^2_{L^2(\hat{\Omega}^s)}
 \end{equation*}
\begin{equation*}
 \le \mathcal{E}_f^n + \mathcal{E}_s^n 
 +  \frac{C}{\mu_f}\Delta t   ||p_{in}(t^n)||^2_{L^2(0,R)}+  \frac{C}{\mu_f}\Delta t   ||p_{out}(t^n)||^2_{L^2(0,R)}.
 \end{equation*}
To get an energy estimate for an arbitrary time-level $\tn \Delta t$, where $1\le \tn\le N$, 
in terms of the energy of the initial data, and the work done by the inlet and outlet dynamic pressure data,
we sum the above inequalities for  $n = 1, \dots,\tn$, cancel the corresponding kinetic energy terms appearing on both sides, and
obtain the energy estimate \eqref{energy_estimate}.
 \end{proof}

\section{Numerical results}\label{sec5} 
We consider two test problems. 
The first was considered in \cite{burman2009stabilization} to test performance of a partitioned scheme 
based on Nitsche's method, with a time penalty term needed for stabilization. The problem involves 
solving a time-dependent Stokes problem coupled with the equations of linear elasticity,
where the coupling is assumed at a fixed fluid domain boundary (linear coupling). 
The second problem we consider was proposed in \cite{quaini2009algorithms} and used in \cite{badiaq2} as a benchmark 
problem for FSI with thick elastic walls to test performance of monolithic FSI schemes. 
In this problem, the full FSI problem  \eqref{sys1}-\eqref{sys2} is solved, and the coupling is evaluated
at the moving fluid-structure interface (nonlinear coupling). 
In both examples, the flow is driven by the time-dependent 
pressure data. To resemble the regime in which instabilities may occure, the fluid and structure densities are taken to be comparable.

We used our operator splitting scheme with $\beta = 1$ to simulate solutions to the two problems.
Our numerical investigations indicate 
that an increase in $\beta\in [0,1]$ increases the accuracy of the scheme. 
Similar experience was also reported in \cite{Martina_paper1} for a FSI problem involving a thin elastic structure. 
This is why all the simulations  presented here correspond to $\beta = 1$.

The value of the viscous parameter $\epsilon$ was taken to be zero in both examples. 
Taking $\epsilon > 0$ regularizes solutions of the FSI problem. 
Therefore,  $\epsilon = 0$ is the most difficult case to consider, since the structure problem in that case is
hyperbolic, exhibiting wave phenomena in the structure at disparate time-scales from the fluid. 

In both examples, the results obtained using the proposed operator splitting scheme are compared with solutions obtained using a monolithic scheme,
showing excellent agreement. In the second example, which considers the full FSI problem  \eqref{sys1}-\eqref{sys2},
we additionally show that our numerical results indicate first-order accuracy in time of the proposed numerical scheme.
Finally,  we show that the condition number of the fluid sub-problem (Problem A2) in the proposed operator splitting scheme 
is by several orders of magnitude smaller than the condition number of monolithic schemes 
due to the fact that no wave  phenomena associated with structural elastodynamics are solved in Problem A2 of the scheme.

\subsection{Example 1. }
We consider a simplified FSI problem in which the fluid is modeled by the time-dependent Stokes problem,
and the structure by the equations of linear elasticity. The fluid-structure coupling is linear in the sense that the fluid 
domain does not change in time. 
The flow is driven by the inlet time-dependent pressure data which is a step function in time:
\begin{equation*}
 p_{in}(t) = \left\{\begin{array}{l@{\ } l} 
10^4 \;\textrm{dyne/cm}^2 & \textrm{if} \; t \le 0.005\\
0 & \textrm{if} \; t>0.005
 \end{array} \right.,   \quad p_{out}(t) = 0 \;\forall t \in (0, T).
\end{equation*} 
The outlet normal stress is kept at zero.

The values of the parameters that determine the fluid and structure geometry,
as well as their physical properties, are given in Table~\ref{T2}
\begin{center}
\begin{table}[ht!]
{\small{
\begin{tabular}{|l l l l |}
\hline
\textbf{Parameters} & \textbf{Values} & \textbf{Parameters} & \textbf{Values}  \\
\hline
\hline
\textbf{Radius} $R$ (cm)  & $0.5$  & \textbf{Length} $L$ (cm) & $5$  \\
\hline
\textbf{Fluid density} $\rho_f$ (g/cm$^3$)& $1.1$ &\textbf{Dyn. viscosity} $\mu$ (poise) & $0.035$    \\
\hline
\textbf{Wall density} $\rho_s $(g/cm$^3$) & $1.2$  & \textbf{Wall thickness} $h$ (cm) & $0.1$  \\
\hline
\textbf{Lam\'e coeff.} $\mu_s $(dyne/cm$^2$) & $5.75 \times 10^5$  & \textbf{Lam\'e coeff.} $\lambda_s $(dyne/cm$^2$) & $1.7 \times 10^6$ \\
\hline
\textbf{Spring coeff.} $\gamma $(dyne/cm$^4$) & $0$  & \textbf{Viscoel. coeff.} $\epsilon $ (dyne s/cm$^2$) & $0$   \\
\hline
\end{tabular}
}}
\caption{Parameter values for Example 1.}
\label{T2}
\end{table}
\end{center}

This problem was suggested in \cite{burman2009stabilization} as a test problem to study performance of an
explicit scheme for FSI problems which was based on Nitsche's method. 
To deal with the instabilities associated with the added mass effect in \cite{burman2009stabilization}, 
a weakly consistent stabilization term was added that
corresponds to the pressure variations at the interface. This decreased the
temporal accuracy of the scheme, which was then corrected by adding certain
defect-correction sub-iterations. 
In \cite{burman2009stabilization} this stabilized explicit scheme was solved using the Taylor-Hood ($\mathbb{P}_2/\mathbb{P}_1$) 
finite elements for the fluid, and $\mathbb{P}_1$ elements for the structure.
The size of the computational mesh was $h_v=0.1$. Their simulations were compared with the solution obtained using a monolithic scheme.
The time step for the monolithic scheme was $\Delta t = 10^{-4}$, while the time step used in the stabilized explicit scheme based on 
Nitsche's method was $\Delta t=10^{-5}$.

In our simulations we discretize the problem in space by using a finite element approach with an isoparametric version of the 
Bercovier-Pironneau element spaces, also known as  the $\mathbb{P}_1$-iso-$ \mathbb{P}_2$ approximation. 
In this approach, a coarse mesh is used for the pressure (mesh size $h_p$), and a fine mesh for velocity (mesh size $h_v=h_p/2$).
To solve the structure problem, we used $\mathbb{P}_1$ elements on a conforming mesh. 
To capture all the physically relevant phenomena, in this example we
chose the space discretization step to be $h_v = \frac{1}{\sqrt{Re}}=0.032$,
where $Re$ is the Reynolds number associated with this problem. 
To achieve comparable accuracy to the scheme proposed in \cite{burman2009stabilization},
we only needed to take the time step $\Delta t = 10^{-4}$, which was the time step used in  \cite{burman2009stabilization} 
for the monolithic solver, but not for the  stabilized explicit scheme proposed there.

\begin{figure}[ht!]
 \centering{
 \includegraphics[scale=0.7]{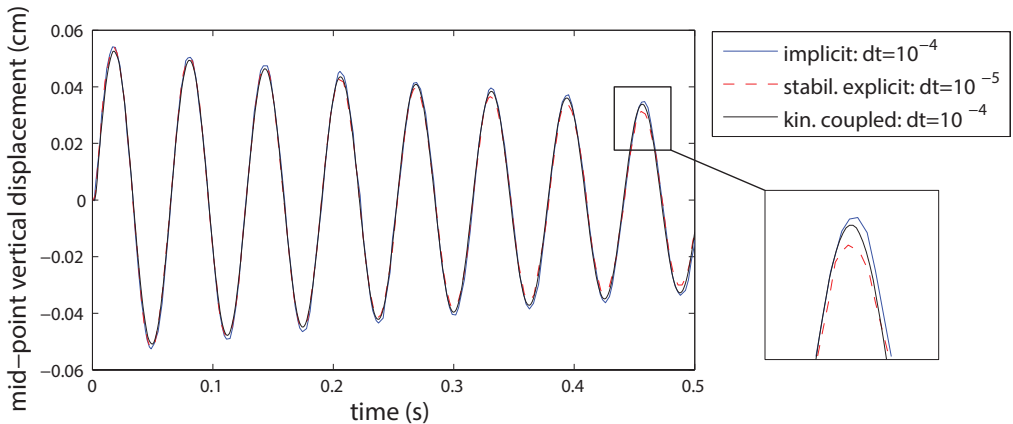}
 }
 \caption{Mid-point vertical displacement computed using a monolithic scheme (blue sold line)
 with $\Delta t = 10^{-4}$, a stabilized explicit scheme proposed by Burman and Fernandez~\cite{burman2009stabilization}
 with $\Delta t = 10^{-5}$ (red dashed line), and the operator splitting scheme with $\beta = 1$ and $\Delta t = 10^{-4}$ (black solid line).}
\label{Ffer1}
 \end{figure}

A comparison between the results of the stabilized explicit scheme proposed in  \cite{burman2009stabilization}, 
the monolithic scheme used in  \cite{burman2009stabilization}, and our operator splitting scheme,
are shown in  Figure~\ref{Ffer1}. In this figure, the displacement of the fluid-structure interface
at the mid-point of the spatial domain was calculated.
Excellent agreement was achieved between our method and the corresponding monolithic method,
where the time step we used was the same as the time step used in the monolithic solver.
In contrast with the stabilized explicit scheme proposed in  \cite{burman2009stabilization}, the operator splitting scheme 
proposed in this manuscript does not require stabilization, providing results of this problem that compare well with the results
of the monolithic scheme using the same time step as in the monolithic solver.

\subsection{Example 2.} 
We consider the fully nonlinear FSI problem  \eqref{sys1}-\eqref{sys2} with the fluid-structure coupling evaluated
at the moving interface (nonlinear coupling). 
The flow is driven by the inlet time-dependent pressure data, which is a cosine pulse lasting for $t_{max}=0.003$ seconds, 
while the outlet normal stress is kept at zero:
 \begin{equation*}
 p_{in}(t) = \left\{\begin{array}{l@{\ } l} 
\frac{p_{max}}{2} \big[ 1-\cos\big( \frac{2 \pi t}{t_{max}}\big)\big] & \textrm{if} \; t \le t_{max}\\
0 & \textrm{if} \; t>t_{max}
 \end{array} \right.,   \quad p_{out}(t) = 0 \;\forall t \in (0, T),
\end{equation*}
where $p_{max} = 1.333 \times 10^4$ (dynes/cm$^2$) and $t_{max} = 0.003$ (s). 
\begin{figure}[ht!]
 \centering{
 \includegraphics[scale=0.60]{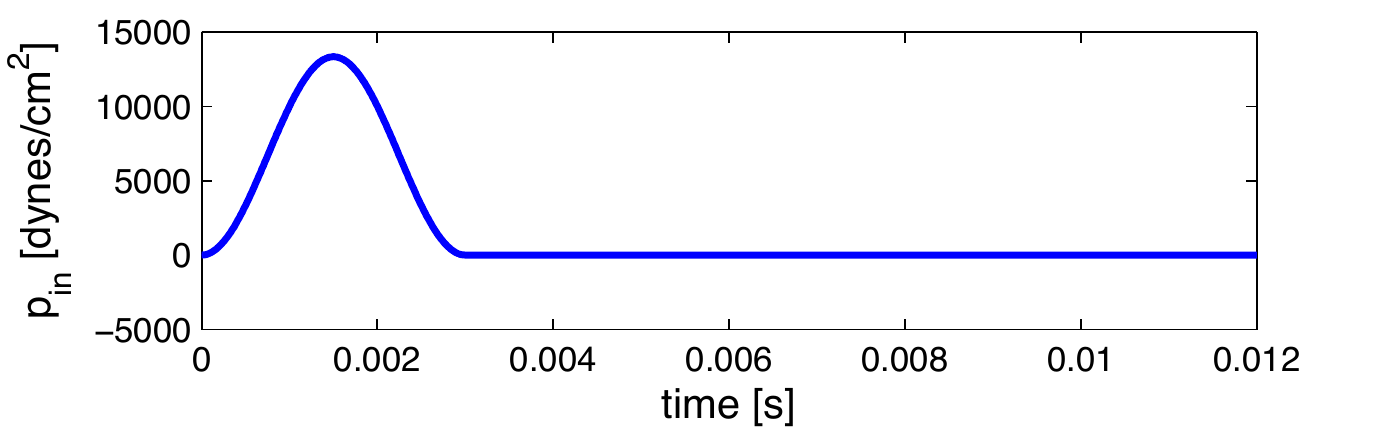}
 }
 \caption{The inlet pressure pulse for Example 2. The outlet normal stress is kept at 0. }
\label{pressure_pulse}
 \end{figure}
 See Figure~\ref{pressure_pulse}. 
 
 This problem was proposed in \cite{quaini2009algorithms} and used in \cite{badiaq2} as a benchmark 
problem for FSI problems in hemodynamics, involving thick elastic walls.

The domain geometry and the values of all the fluid and structure parameters in this example are given in Table~\ref{T1}.
{{
\begin{center}
\begin{table}[ht!]
{\small{
\begin{tabular}{|l l l l |}
\hline
\textbf{Parameters} & \textbf{Values} & \textbf{Parameters} & \textbf{Values}  \\
\hline
\hline
\textbf{Radius} $R$ (cm)  & $0.5$  & \textbf{Length} $L$ (cm) & $6$  \\
\hline
\textbf{Fluid density} $\rho_f$ (g/cm$^3$)& $1$ &\textbf{Dyn. viscosity} $\mu$ (g/cm s) & $0.035$    \\
\hline
\textbf{Wall density} $\rho_s $(g/cm$^3$) & $1.1$  & \textbf{Wall thickness} $h$ (cm) & $0.1$  \\
\hline
\textbf{Lam\'e coeff.} $\mu_s $(dyne/cm$^2$) & $5.75 \times 10^5$  & \textbf{Lam\'e coeff.} $\lambda_s $(dyne/cm$^2$) & $1.7 \times 10^6$ \\
\hline
\textbf{Spring coeff.} $\gamma $(dyne/cm$^4$) & $4 \times 10^6$  & \textbf{Viscoel. coeff.} $\epsilon $ (dyne s/cm$^2$) & $0$   \\
\hline
\end{tabular}
}}
\caption{Parameter values for Example 2.}
\label{T1}
\end{table}
\end{center}
}}

Problem \eqref{sys1}-\eqref{sys2} was solved over the time interval $[0,0.012]$s.
This is the time that it takes the pressure pulse generated at the inlet,
to travel across the entire fluid domain and reach the outlet. 

We solved this problem using the following 5 different time steps: $\Delta t = 10^{-4}, 5 \times 10^{-5}, 10^{-5}$, $5 \times 10^{-6}$,
and $10^{-6}$.
These time steps were used to numerically show that the method is convergent, and that its accuracy in time is first-order.

 \begin{figure}[ht!]
\center{
\includegraphics[scale=0.5]{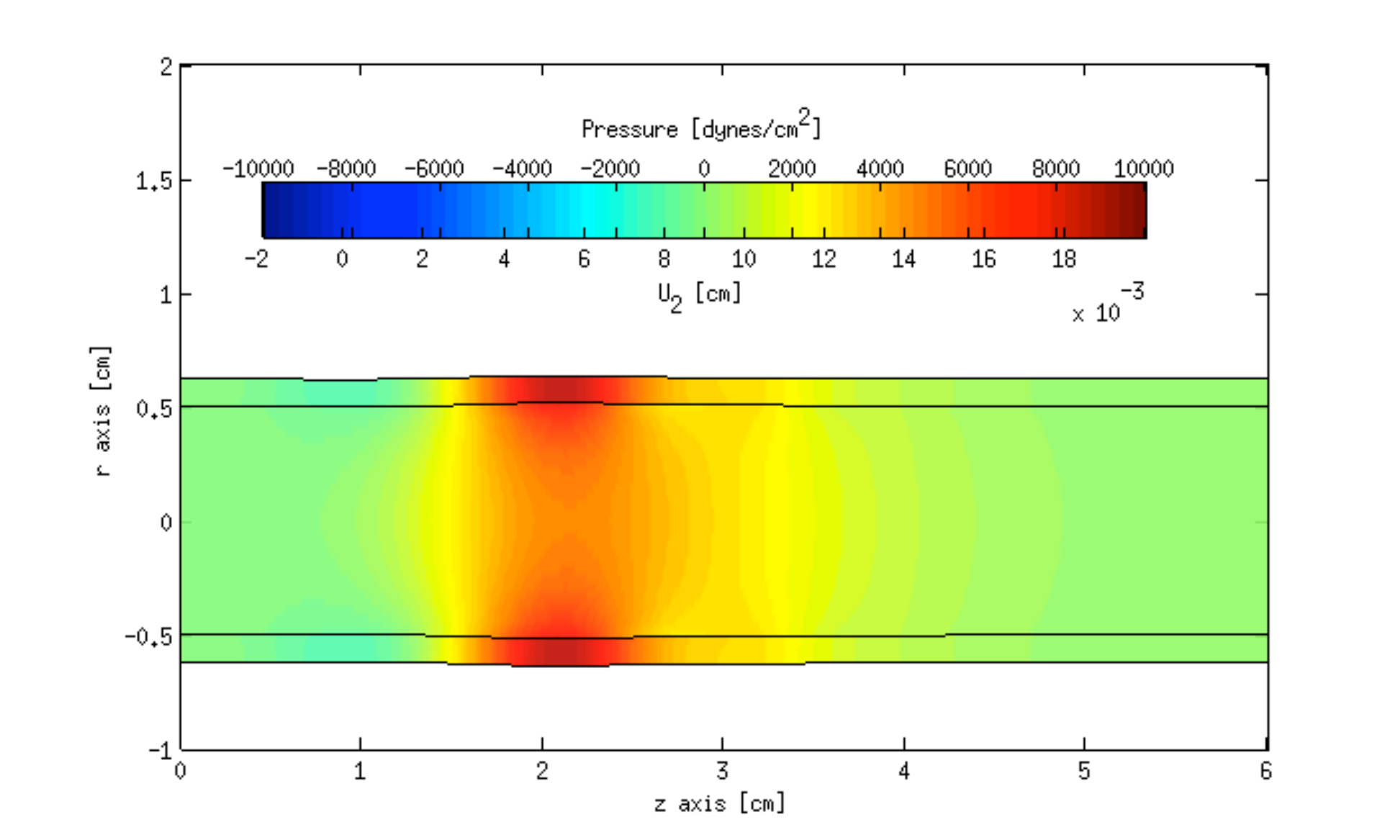}
}
 \caption{A snapshot of the pressure wave traveling from left to right coupled with the radial component of the structure displacement.
 The legend shows the values for the pressure (top scale) and displacement (bottom scale) over the same color scale. }
\label{2Dpressure}
 \end{figure}
 
  \begin{figure}[ht!]
\center{
\includegraphics[scale=0.5]{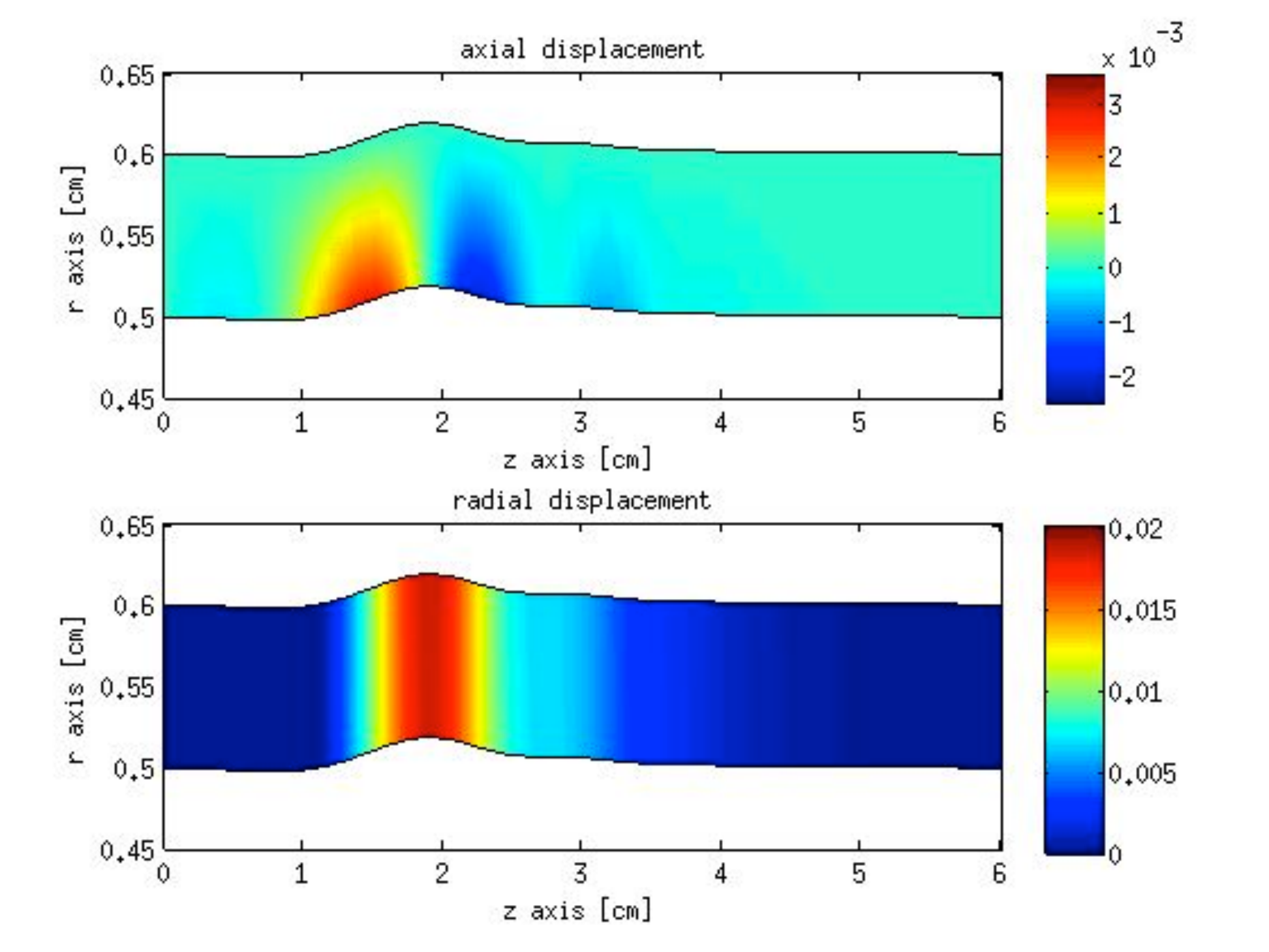}
}
 \caption{The longitudinal (top) and radial (bottom) components of the structure displacement. Notice how the red and blue colors in 
 longitudinal displacement denote longitudinal stretching of the structure in opposite (+/-) directions. }
\label{2Ddisplacement}
 \end{figure}
 
 As before, to discretize the problem in space, we used the 
Bercovier-Pironneau element spaces, also known as  the $\mathbb{P}_1$-iso-$ \mathbb{P}_2$ approximation, with 
the velocity mesh of size  $h_v = 0.01$, and the pressure mesh of size $h_p=2 h_v=0.02$.
To solve the structure problem, the $\mathbb{P}_1$ elements were used on a conforming mesh. 

Figures~\ref{2Dpressure} and \ref{2Ddisplacement} show 2D plots of the fluid pressure and structure displacement, respectively. 
The pressure wave, shown in Figure~\ref{2Dpressure}, travels from left to right, displacing the thick structure. The colors of the
thick structure displacement in Figure~\ref{2Dpressure} denote the magnitude of the radial component of displacement. 
Figure~\ref{2Ddisplacement} shows separate snapshots of the longitudinal and radial components of the displacement.
Notice how the red and blue colors in longitudinal displacement denote longitudinal stretching of the structure in opposite directions.

\begin{figure}[ht!]
 \centering{
 \includegraphics[scale=0.70]{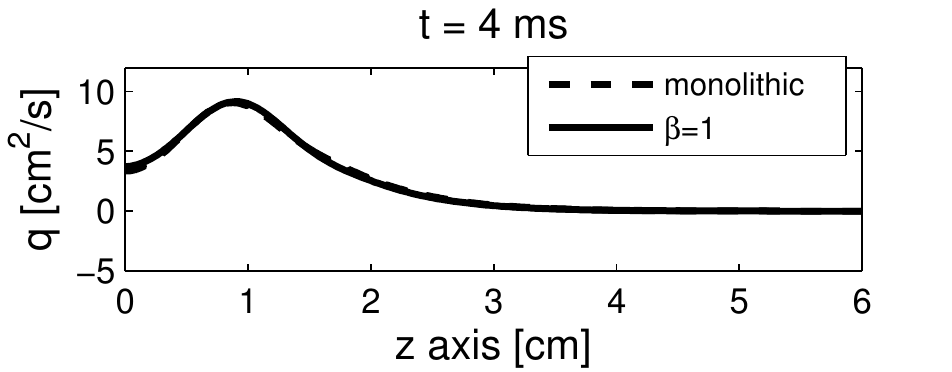}
\includegraphics[scale=0.70]{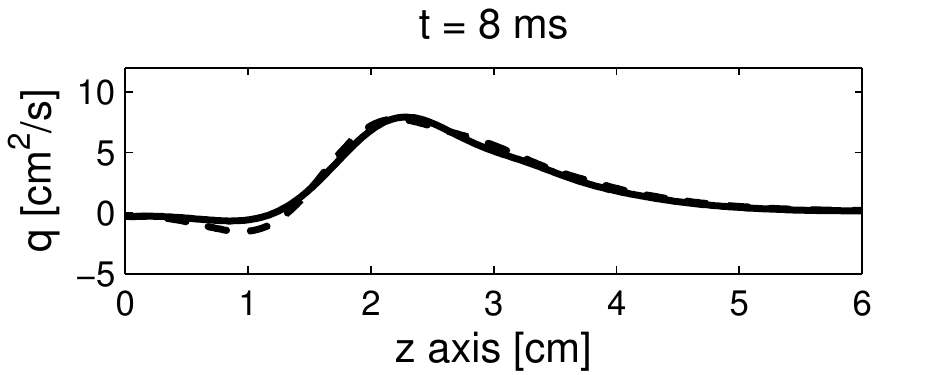}
\includegraphics[scale=0.70]{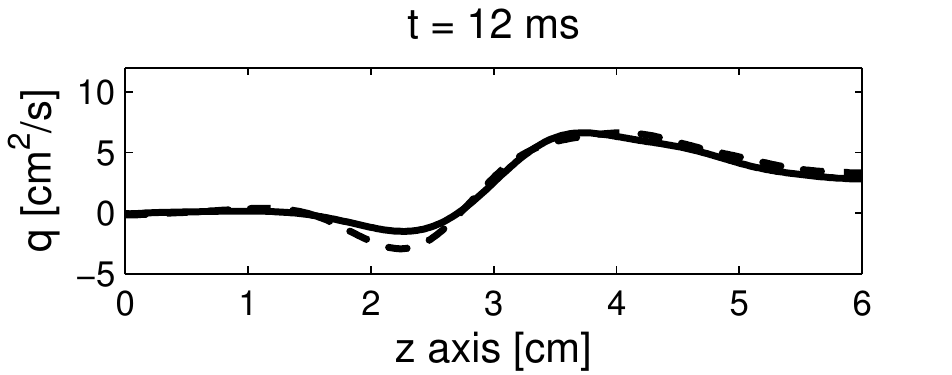}

 }
 \caption{Flow rate vs. z, at $t=4,8,12$ ms,
  computed with the monolithic scheme by Quaini~\cite{quaini2009algorithms} (time step $\Delta t = 10^{-4}$; dashed red line), and 
 with our operator splitting scheme with $\beta = 1$ (time step $\Delta t =5 \times 10^{-5}$; solid blue line). }
\label{flow}
 \end{figure}
\begin{figure}[ht!]
 \centering{
    \includegraphics[scale=0.70]{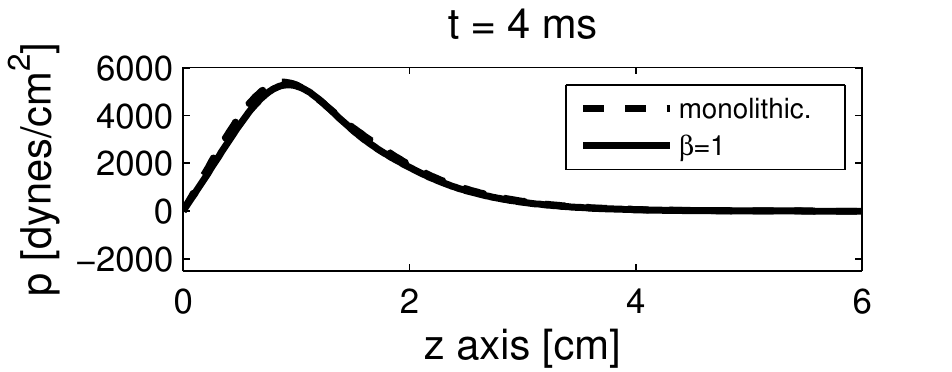}
\includegraphics[scale=0.70]{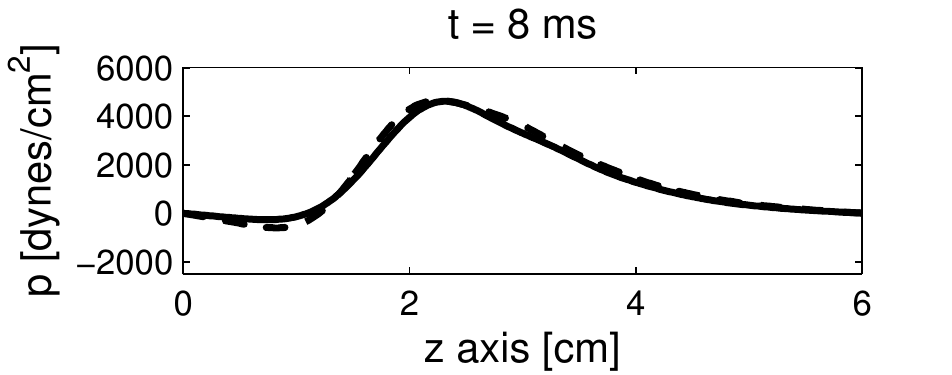}
\includegraphics[scale=0.70]{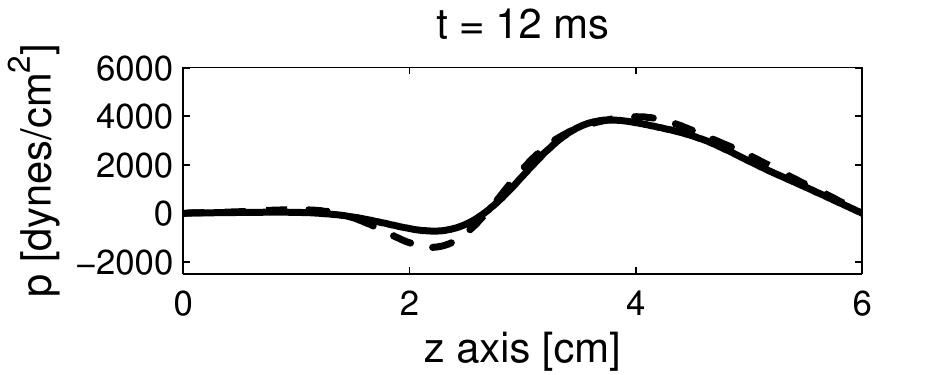}

 }
 \caption{ Mean pressure vs. z, at $t=4,8,12$ ms,
  computed with the monolithic scheme by Quaini~\cite{quaini2009algorithms} (time step $\Delta t = 10^{-4}$; dashed red line), and 
 with our operator splitting scheme with $\beta = 1$ (time step $\Delta t =5 \times 10^{-5}$; solid blue line). }
\label{press}
 \end{figure}
 \begin{figure}[ht!]
 \centering{
\includegraphics[scale=0.70]{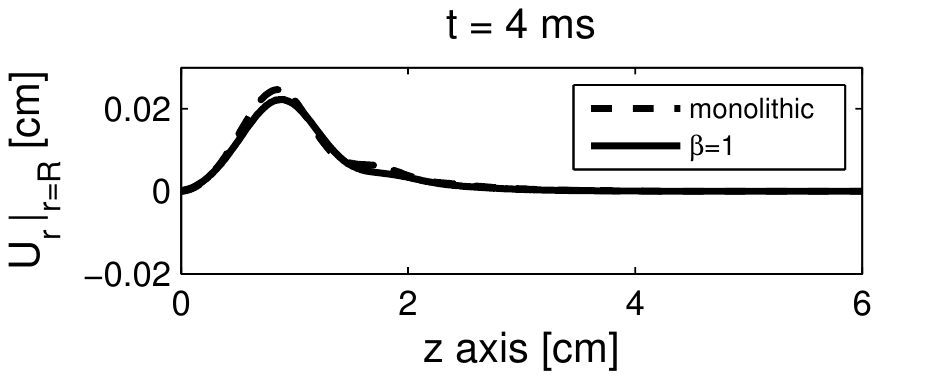}
\includegraphics[scale=0.70]{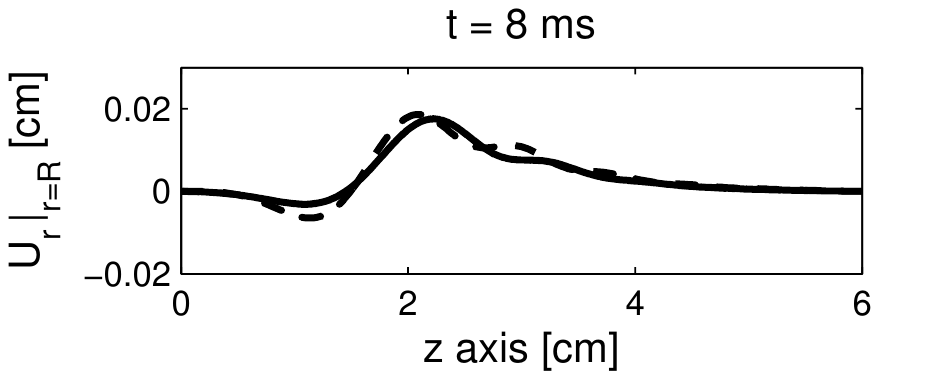}
\includegraphics[scale=0.70]{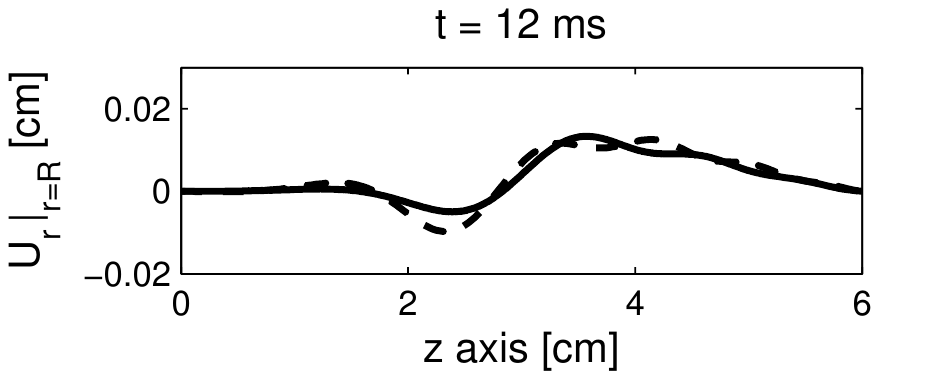}
 }
 \caption{Fluid-structure interface displacement vs. z, at $t=4,8,12$ ms,
  computed with the monolithic scheme by Quaini~\cite{quaini2009algorithms} (time step $\Delta t = 10^{-4}$; dashed red line), and 
 with our operator splitting scheme with $\beta = 1$ (time step $\Delta t =5 \times 10^{-5}$; solid blue line). }
\label{displfig}
 \end{figure}

The numerical results obtained using the operator splitting scheme proposed in this manuscript
were compared with the numerical results obtained using the monolithic scheme proposed in~\cite{quaini2009algorithms,badiaq2}.
The monolithic scheme was solved on the same mesh using stabilized $\mathbb{P}_1-\mathbb{P}_1$ elements for the fluid problem, and 
$\mathbb{P}_1$ elements for the structure problem. 

Figures~\ref{flow},~\ref{press}, and~\ref{displfig}, show the calculated
flow rate, mean pressure, and fluid-structure interface displacement, respectively, 
as functions of the horizontal axis $z$, for three different snap-shots.
The three figures show a comparison between our operator splitting scheme, shown in solid line, 
and the monolithic scheme of \cite{quaini2009algorithms,badiaq2}, shown in dashed line. 
In these figures, the time step used in the monolithic scheme was $\Delta t = 10^{-4}$. To obtain roughly the same
accuracy, we used the time step $\Delta t / 2 = 5 \times 10^{-5}$. 
Figures~\ref{flow},~\ref{press}, and~\ref{displfig}  show a good comparison between the two solutions.

\subsubsection{Convergence in time}\label{accuracy}
\begin{figure}[ht!]
 \centering{
 \includegraphics[scale=0.77]{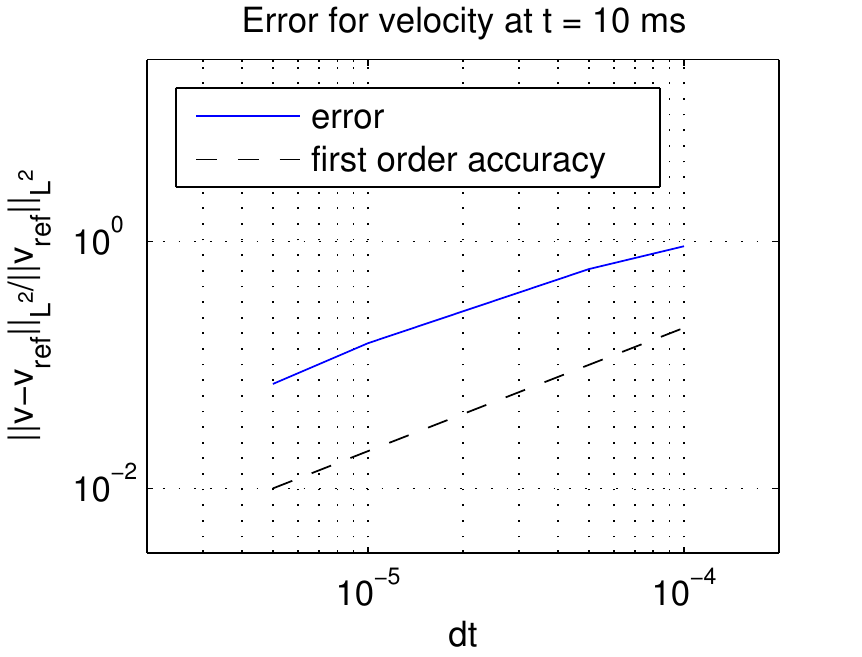}
\includegraphics[scale=0.77]{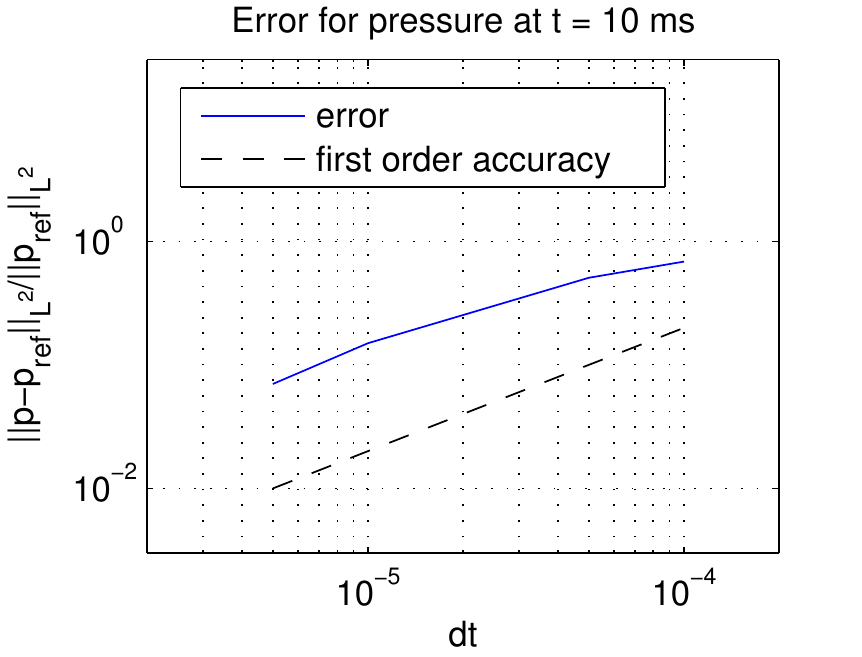} 
\includegraphics[scale=0.77]{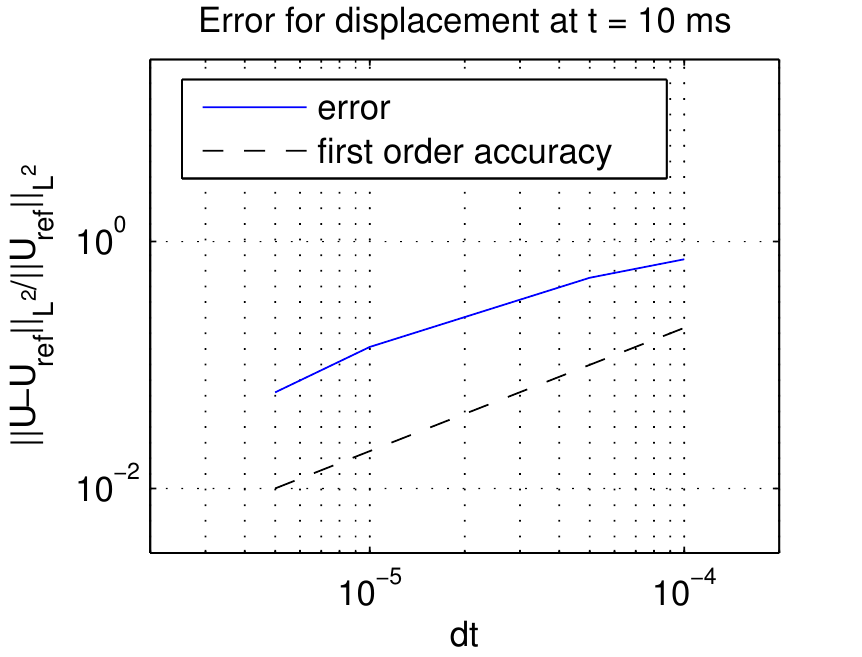}
 }
 \caption{Example 1: Figures show relative errors obtained at t=10 ms. Top left: Relative error for fluid velocity. Top right: Relative error for fluid pressure. Bottom: Relative error for the structure displacement.}
\label{error}
 \end{figure}

In this example, we also study convergence in time of our operator-splitting scheme. 
For this purpose, we define the reference solution to be the one obtained with $\Delta t = 10^{-6}$. 
We calculated the relative $L^2$-errors for the velocity, pressure and displacement,
between the reference solution and the solutions obtained using 
$\Delta t = 10^{-4}, 5 \times 10^{-5},\;  10^{-5}$,  and $5 \times 10^{-6}$.
Figure~\ref{error} shows the log-log plot of the relative errors, superimposed over a line with slope $1$, corresponding to first-order accuracy.
The slopes indicate that our scheme is first-order accurate in time. Indeed, 
Table~\ref{T4} shows the precise numbers from Figure~\ref{error}, calculated at time $t = 10$ ms,
indicating first-order in time convergence for the velocity, pressure, and displacement.

 \begin{table}[ht!]
\begin{center}
{\small{
\begin{tabular}{| l  c  c  c  c  c  c |}
\hline
$ \Delta t $ & $\frac{||p-p_{ref}||_{L^2}}{||p_{ref}||_{L^2}} $ & $L^2$ order & $\frac{||\boldsymbol u-\boldsymbol u_{ref}||_{L^2}}{||\boldsymbol u_{ref}||_{L^2}}$  &$L^2$ order & $ \frac{||\boldsymbol U - \boldsymbol U_{ref}||_{L^2}}{|| \boldsymbol U_{ref}||_{L^2}} $ & $L^2$ order \\
\hline
\hline
$10^{-4}$ & $0.69$  & - & $0.92$  & - & $0.72$  & - \\ 
\hline  
$5 \times 10^{-5}$ & $0.51$  & 0.43 & $0.60$  & 0.60 & $0.51$  & 0.5 \\ 
\hline   
$10^{-5}$ & $0.15$  & 0.76 & $ 0.15$  & 0.85 & $0.14$  & 0.8 \\ 
\hline   
$5 \times 10^{-6}$ & $0.07$  & 1.06 & $0.07$  & 1.1 & $0.06$  & 1.1 \\ 
\hline 
\end{tabular}
}}
\end{center}
\caption{Convergence in time calculated at $t = 10$ ms. }
\label{T4}
\end{table}

\subsubsection{The condition number of the fluid sub-problem in Problem A2}\label{condition_number}
We recall that in Problem A2 of this operator splitting scheme a fluid sub-problem is solved in such a way to include the structure inertia
into the fluid sub-problem. This was done for stability reasons, i.e., to avoid issues related to the added mass effect, associated 
with classical Dirichlet-Neumann partitioned schemes for FSI in hemodynamics.
As mentioned earlier, 
for FSI problems containing a thin fluid-structure interface with mass, including the fluid-structure interface inertia into the fluid sub-problem
can be easily accomplished
via a Robin-type boundary condition for the fluid sub-problem, leading to a fully partitioned scheme \cite{Martina_paper1,BorSunMulti,Martina_Multilayered}. In our problem, however, the fluid-structure interface is just a trace of the thick structure in contact with the fluid.
In this case, the trace of the fluid-structure interface does not have a well-defined inertia, and we need to include the inertia of the entire
thick structure into the fluid sub-problem. As we saw earlier, this is done by solving a fluid sub-problem together with a simple 
problem for the structure that takes only the structure inertia into account (and possibly the viscous effects of the structure if $\epsilon \ne 0$).
Thus, in Problem A2, we solve a simplified coupled problem consisting of a fluid sub-problem, and a structure sub-problem involving only
structure inertia, coupled through a simple continuity of stresses condition. Even though this is reminiscent of monolithic schemes, 
we show here that the condition number of this fluid sub-problem is by several orders of magnitude smaller than the condition number
of the monolithic scheme by Quaini \cite{quaini2009algorithms,badiaq2}.
Indeend, we calculated the condition number $\kappa$ of the stiffness matrix  in Problem A2 of our operator splitting scheme and obtained that
$\kappa_{\rm split} = 1.08\times 10^4$. Similarly, we calculated
the related condition number for the stiffness matrix of the monolithic scheme by Quaini~\cite{quaini2009algorithms,badiaq2}
and obtained that it equals $\kappa_{\rm mono} = 7.82\times 10^8$.
This is directly related to the fact that in Problem A2 of our scheme, no structural elastodynamics problem was solved that would capture
the wave phenomena in the structure traveling at disparate time-scales form the pressure wave in the fluid.

Therefore, although the problem in Problem A2 is solved
on both fluid and structure domains, its condition number is equivalent to that of pure fluid solvers.
Thus, the proposed operator-splitting scheme consists of a fluid module and a structure module, which can be easily replaced
if different structure models or different solvers are to be used, see, e.g., \cite{Lukacova}.
Furthermore, for more general multi-physics problems, 
additional modules can be easily added to capture different physics in the problem, as was done
in \cite{Martina_Multilayered} to study FSI with multi-layered structures, in \cite{Martina_Biot} to study FSI with multi-layered
poroelastic walls, or in \cite{BorSunStent} to include a model of a stent in the underlying FSI problem.
Modularity, unconditional stability, and simple implementation are the features that make this 
operator-splitting scheme particularly appealing for multi-physics problems involving fluid-structure interaction.

\if 1 = 0
\subsection{Dependence on $\epsilon$}
In the previous example we compared our numerical results with results obtained by a monolithic scheme in which structure is assumed to be purely elastic. With our visicolestic parameter $\epsilon=0.01$, we obtained very good agreement with a monolithic scheme. In this subsection we investigate dependance on viscoelastic parameter $\epsilon$. 

Figures~\ref{Fig_gamma_vel},~\ref{Fig_gamma_press}, and~\ref{Fig_gamma_displ} show the comparison of results obtained with our scheme using values $\epsilon=1, \epsilon=0.1$, and $\epsilon=0.01$.
\begin{figure}[ht!]
 \centering{
\includegraphics[scale=0.8]{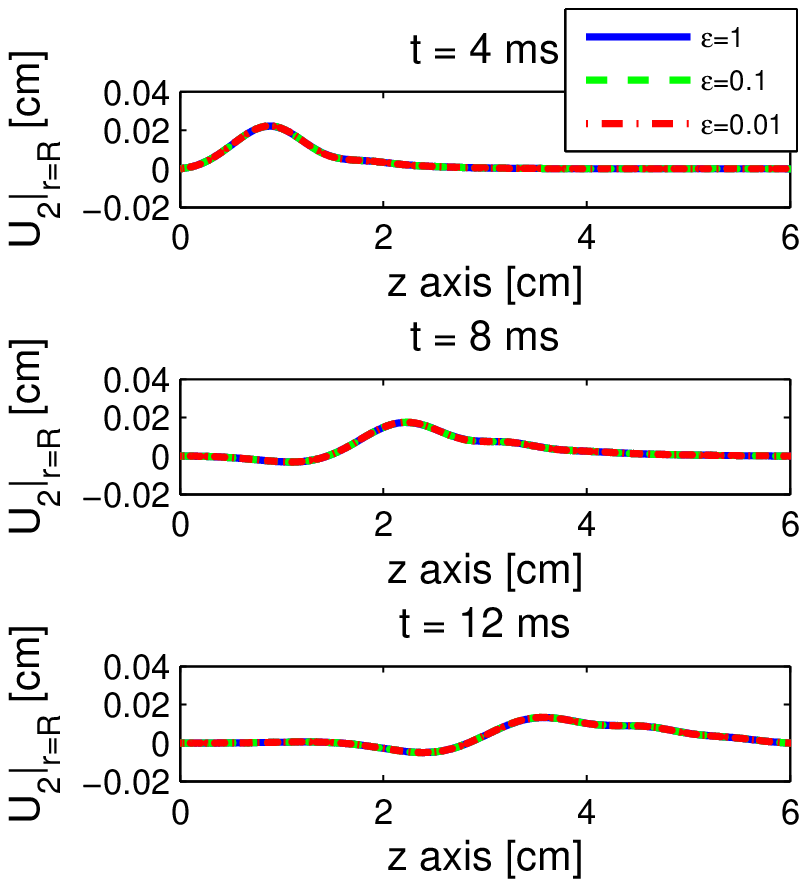}
 }
 \caption{ Inner radius of the tube computed using values $\epsilon=1$ (solid blue line), $\epsilon=0.1$  (dashed green line), and $\epsilon=0.01$ (dash-dot red line). }
\label{Fig_gamma_displ}
 \end{figure}
\begin{figure}[ht!]
 \centering{
 \includegraphics[scale=0.8]{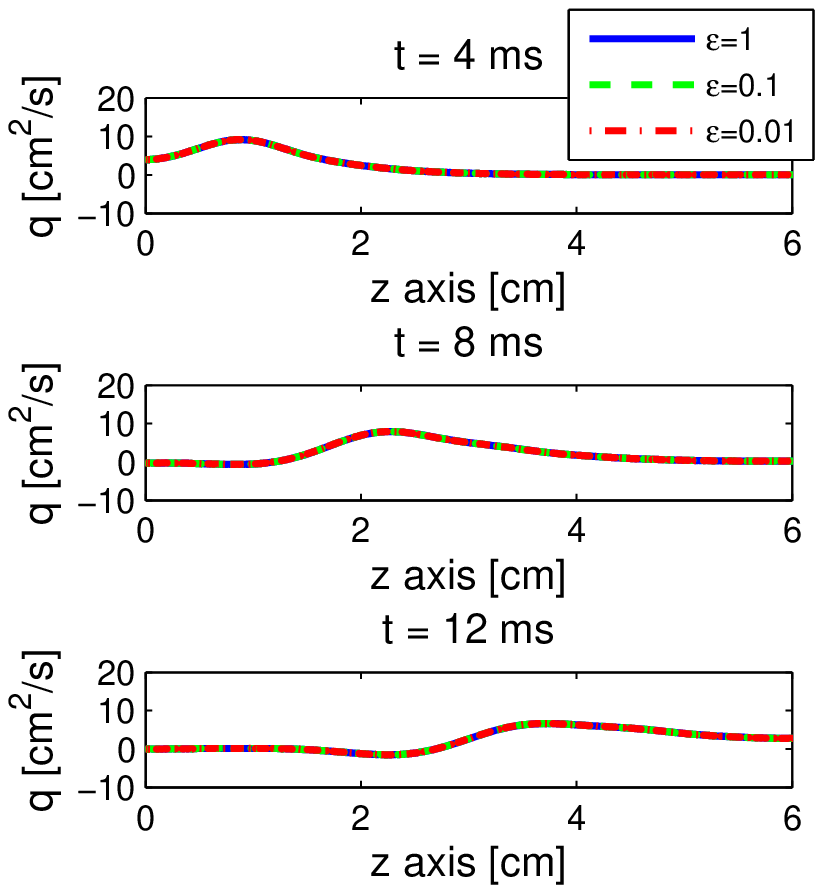}
 }
 \caption{Flowrate computed using values $\epsilon=1$ (solid blue line), $\epsilon=0.1$  (dashed green line), and $\epsilon=0.01$ (dash-dot red line). }
\label{Fig_gamma_vel}
 \end{figure}
\begin{figure}[ht!]
 \centering{
 \includegraphics[scale=0.8]{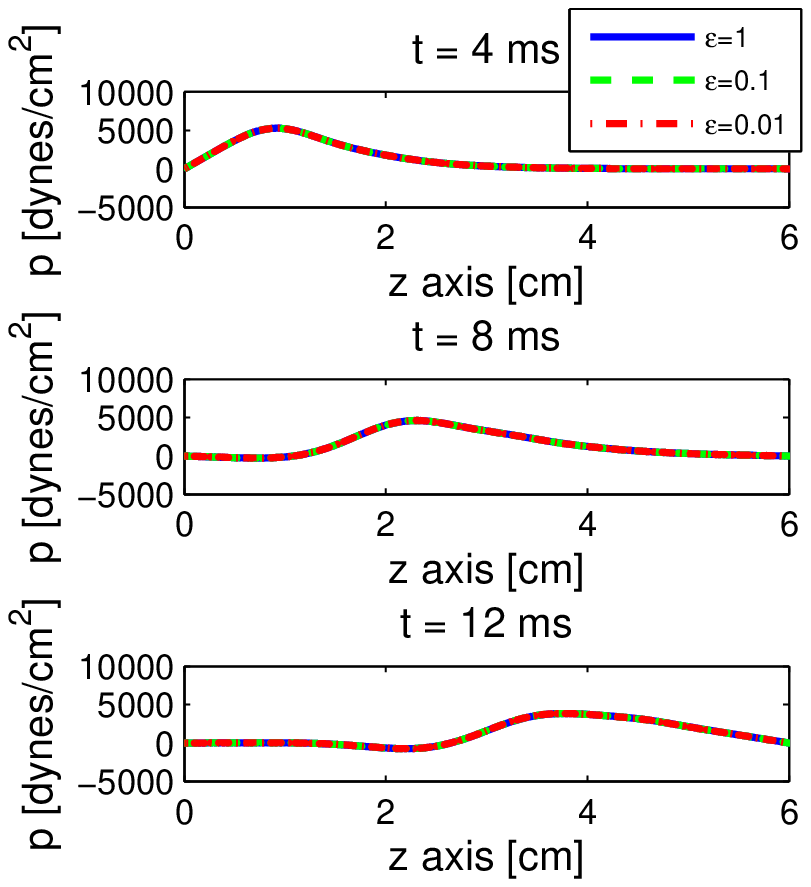}
 }
 \caption{Mean pressure computed using values $\epsilon=1$ (solid blue line), $\epsilon=0.1$  (dashed green line), and $\epsilon=0.01$ (dash-dot red line). }
\label{Fig_gamma_press}
 \end{figure}

Taking the results obtained with $\epsilon=0.01$ as a reference solution,  relative errors for velocity, pressure and displacement at t=12 ms, between the reference solution and solutions obtained with $\epsilon =1$ and $\epsilon=0.1$, are given in Table~\ref{Table_error_gamma}.
 \begin{table}[ht!]
\begin{center}
{\scriptsize{
\begin{tabular}{| l  c  c  |}
\hline
$  $ & $\epsilon=1 $ &  $\epsilon=0.1$    \\
\hline
\hline
$\displaystyle\frac{||p-p_{ref}||_{L^2}}{||p_{ref}||_{L^2}}$ & $0.0064$ & $ 6.28 \times 10^{-4}$   \\ 
\hline  
$\displaystyle\frac{||\boldsymbol u-\boldsymbol u_{ref}||_{L^2}}{||\boldsymbol u_{ref}||_{L^2}}$ & $0.007$ & $ 6.78 \times 10^{-4}$   \\ 
\hline   
$\displaystyle\frac{||\boldsymbol \eta - \boldsymbol \eta_{ref}||_{L^2}}{||\boldsymbol \eta_{ref}||_{L^2}}$  & 0.0083 & $7.97 \times 10^{-4}$\\ 
\hline 
\end{tabular}
}}
\end{center}
\caption{Relative errors for pressure, velocity, and displacement between solution obtained with $\epsilon=0.01$ (reference solution), and solutions  obtained with $\epsilon =1$ and $\epsilon=0.1$.}
\label{Table_error_gamma}
\end{table}
\fi


\section{Conclusions}
This work proposes a modular scheme for fluid-structure interaction problems with thick structures. The proposed scheme is
based on the Lie operator splitting approach, which separates the fluid from the structure sub-problem, and on using an Arbitrary
Lagrangian Eulerian approach to deal with the motion of the fluid domain. To achieve unconditional stability without 
sub-iterations in each time step,
the fluid sub-problem includes structure inertia, which requires solving the fluid sub-problem on both domains, i.e., the fluid and structure domains. 
While this is reminiscent of monolithic schemes,  the condition number of the fluid sub-step is significantly smaller than the condition number 
associated with classical 
FSI  monolithic schemes. This is because the wave propagation in the elastic structure is treated separately in the structure sub-problem,
and not together with the fluid problem, as in classical monolithic schemes.
The advantage of this approach over classical monolithic schemes is the possibility to use larger time steps, the separation of 
dissipative vs. non-dissipative features of the coupled problem allowing the use of non-dissipative
solvers to treat wave propagation in the structure, and modularity, which allows simple extensions of the scheme to capture 
different multi-physics problems associated with FSI. A disadvantage of this scheme over classical partitioned schemes is that 
the fluid sub-step requires solving the associated problem on both the fluid and structure domains in a monolithic fashion. 
However, unlike the classical
Dirichlet-Neumann partitioned schemes, the proposed scheme is unconditionally stable for all the parameters in the problem.
This was shown in the present manuscript by proving an energy estimate associated with unconditional stability of the scheme, 
for the full, nonlinear FSI problem. 

 \section{Acknowledgements} The authors would like to thank the reviewers for their careful reading of the manuscript and for the 
 insightful suggestions that improved the quality of the manuscript.
 
 This research has been supported in part by the National Science Foundation under grants
DMS-1318763 (\v{C}ani\'{c} and Buka\v{c}),
 DMS-1311709 (\v{C}ani\'{c} and Muha), 
 DMS-1262385 and DMS-1109189 (\v{C}ani\'{c} and Quaini), 
 DMS-0806941 (\v{C}ani\'{c} and Buka\v{c}),
 and by the Texas Higher Education Board (ARP-Mathematics) 003652-0023-2009 (\v{C}ani\'{c} and Glowinski).

\bibliographystyle{plain}
\bibliography{thick}
\end{document}